\documentclass[11pt]{article}

\usepackage[letterpaper,
            left=1in,
            textwidth=6.6in,      % right margin becomes 2.4in
            textheight=8.6in,
            vmarginratio=1:1,
            headheight=14pt,
            footskip=0.4in]{geometry}

\usepackage[T1]{fontenc}
\usepackage[utf8]{inputenc}
\usepackage[largesc]{newpxtext}          % Palatino-like text, real small caps
\usepackage[vvarbb]{newpxmath}           % matching math, blackboard bold
\usepackage[scaled=0.93]{inconsolata}    % code, if any
\usepackage{microtype}                   % better spacing; no visible cost
\definecolor{UCIblue}{gray}{0.5}
\usepackage{amsmath,amssymb,amsthm,mathtools}
\usepackage{booktabs,array}
\usepackage[inline]{enumitem}
\setlist{itemsep=0.2\baselineskip,topsep=0.3\baselineskip,leftmargin=2em}
\usepackage{graphicx,tikz}
\usetikzlibrary{patterns,decorations.pathreplacing}

\usepackage{xcolor}
\usepackage{caption}
\usepackage{subcaption}
\definecolor{accent}{RGB}{47,79,118}     % muted slate blue
\definecolor{accentlight}{RGB}{237,241,246}
\definecolor{rule}{RGB}{160,170,185}

\usepackage{titlesec}
\titleformat{\section}
  {\Large\bfseries\color{accent}}{\thesection.}{0.6em}{}
\titleformat{\subsection}
  {\large\bfseries\color{accent}}{\thesubsection.}{0.6em}{}
\titleformat{\subsubsection}[runin]
  {\bfseries}{\thesubsubsection.}{0.5em}{}[.]
\titlespacing*{\section}{0pt}{1.6\baselineskip}{0.7\baselineskip}
\titlespacing*{\subsection}{0pt}{1.2\baselineskip}{0.5\baselineskip}
\titlespacing*{\subsubsection}{0pt}{0.8\baselineskip}{0.6em}

\titleformat{\paragraph}[runin]{\bfseries}{}{0pt}{}[.]
\titlespacing*{\paragraph}{0pt}{0.6\baselineskip}{0.6em}

\theoremstyle{plain}

\theoremstyle{definition}

\theoremstyle{remark}

\usepackage[most]{tcolorbox}
\tcbset{
  tutorialbox/.style={
    enhanced, breakable,
    colback=accentlight, colframe=accent,
    boxrule=0.5pt, arc=1.5pt,
    left=8pt, right=8pt, top=6pt, bottom=6pt,
    fonttitle=\bfseries\small, coltitle=white,
    colbacktitle=accent,
    attach boxed title to top left={yshift=-2mm, xshift=4mm},
    boxed title style={arc=1pt, boxrule=0pt, left=5pt, right=5pt, top=2pt, bottom=2pt},
    before skip=0.9\baselineskip, after skip=0.9\baselineskip,
  }
}
\newtcolorbox{algobox}[1]{tutorialbox, title={#1}}
\newtcolorbox{summarybox}[1]{tutorialbox, title={#1}}

\numberwithin{equation}{section}
\numberwithin{table}{section}
\numberwithin{figure}{section}

\usepackage[font=small,labelfont=bf,labelsep=period,width=0.92\textwidth]{caption}

\usepackage[colorlinks=true,
            linkcolor=accent, citecolor=accent, urlcolor=accent]{hyperref}
\hypersetup{
  pdftitle={From Inverse Problems to Adjoint Gradients: A Discrete, Implementation-Oriented Tutorial on PDE-Constrained Optimization},
  pdfauthor={Solmaz Kia},
  colorlinks=true,
  linkcolor=accent, citecolor=accent, urlcolor=accent
}
\usepackage{fancyhdr}
\fancypagestyle{plain}{\fancyhf{}\fancyfoot[C]{\small\thepage}
  }

\usepackage{authblk}

\title{\bfseries\color{accent}
  From Inverse Problems to Adjoint Gradients:\\[0.25em]
  \large\mdseries\color{black}
  A Discrete, Implementation-Oriented Tutorial on\\
  PDE-Constrained Optimization}
\author{Solmaz Kia}
\affil{Department of Mechanical and Aerospace Engineering,
       University of California, Irvine\\
       \texttt{solmaz@uci.edu}}
\date{}

\renewenvironment{abstract}
  {\begin{center}\begin{minipage}{0.88\textwidth}\small
   \noindent\textbf{Abstract.}\ }
  {\end{minipage}\end{center}\vspace{0.5\baselineskip}}

\usepackage{tocloft}
\usepackage{setspace}
\begin{document}
\vspace{-1.5in}
\maketitle
\thispagestyle{plain}

\begin{abstract}
The adjoint method is a standard approach for computing gradients in
PDE-constrained optimization and central to parameter estimation, model calibration, and optimal
control.
This tutorial presents a discrete, implementation-oriented introduction
to the method with an emphasis on conceptual understanding. Starting
from an inverse problem governed by a partial differential equation, it
develops the corresponding optimization problem and derives the adjoint
gradient from the perspective of iterative descent methods. The roles of
the state, the parameters, the objective, and the constraint are made
explicit, and the adjoint equation is obtained directly from the
gradient computation rather than postulated. Particular attention is
given to explicit derivations, implementation details, and the
connection between the mathematical development and the resulting
computational algorithm. The tutorial also examines the influence of
discretization on the adjoint construction, explains the origin of the
backward-in-time adjoint recursion, and presents practical procedures
for gradient verification. A vibrating-string inverse problem serves as
a running example, illustrating the complete workflow from model
formulation to adjoint-based optimization. The presentation is intended
for advanced undergraduate students, graduate students, and
practitioners seeking a step-by-step understanding of adjoint gradient~computation.
\end{abstract}
\vspace{-0.3in}
\begingroup
\setstretch{0.9}
\tableofcontents
\endgroup
\section{Introduction}

The adjoint method is the standard tool for computing gradients in
optimization problems governed by partial differential equations, and it is
central to parameter estimation, model calibration, and optimal
control~\cite{lions,Giles2000,Borzi2012}. It computes the gradient of an
objective function with respect to all parameters at a cost comparable to a
single simulation, largely independent of the number of parameters, and it is
this property that makes large-scale gradient-based optimization practical.
The method is most naturally understood through the relationship between two
problems: the \emph{forward problem}, which determines the system state from
prescribed data, and the \emph{inverse problem}, which seeks unknown
parameters from observations. This tutorial develops the connection between
these two problems and traces the construction of the adjoint method from the
inverse problem formulation to a computable~gradient.

In the forward problem, the governing equations are solved for the state of
the system given its material properties, geometry, sources, and boundary
data. This is the problem addressed by numerical simulation codes, and because
the equations must be resolved on a discretization fine enough to capture the
relevant physics, even a single solve may be expensive. The inverse problem
runs in the opposite direction, and it is the direction in which most problems
of practical interest are posed: some of the quantities that specify the
system---a material property distributed through the domain, the strength or
location of a source, a boundary condition, the geometry, a control
input---are not known and must be inferred from measurements or chosen by
design. These are the \emph{parameters}. The map from parameters to
observations is typically nonlinear, unstable, and not available in closed
form~\cite{Borzi2012,Tarantola2005,plessix}, so the inverse problem is posed
instead as an optimization problem: an objective measures the discrepancy
between observations and the states predicted by a candidate set of
parameters, and it is minimized subject to the governing equations. The
complication that governs everything that follows is that the objective does
not depend on the parameters directly, but through the solution of the
governing equations; the effect of a parameter perturbation on the objective
can therefore be determined only through the corresponding change in the
forward solution.

This constrained minimization problem is solved iteratively, and the choice of method is dictated by dimension. Methods that use objective values
alone~\cite{conn2009,larson2019}, including the ensemble-based variants used
for inverse problems~\cite{iglesias2013}, require nothing more than the
ability to run the forward simulation, but the number of evaluations they need
grows rapidly with the number of parameters, which rules them out once the
unknown is a discretized field. Large-scale optimization therefore relies on first-order methods, which require the gradient of the objective with respect to the parameters. Because of the indirect dependence just described, that gradient is not available in closed form, and approximating it by finite differences requires one forward solve per parameter~\cite{nocedal,plessix}---prohibitive when the unknown carries thousands to millions of degrees of freedom. The adjoint method removes this
dependence on the number of parameters: by introducing an auxiliary variable,
the gradient is obtained from the solution of a single adjoint problem whose cost is comparable to that of the forward simulation.\footnote{In the number of simulations required, not in memory; see Section~\ref{sec:backward}'s remark on storage and checkpointing.}

The consequences of this single idea have been substantial. In aerodynamics,
Jameson's formulation of design as a control problem, with the gradient
supplied by an adjoint flow solver, made it possible to optimize wing shapes
described by hundreds of design variables with each gradient costing roughly
one flow solution and one adjoint solution, regardless of the number of
variables, and it remains the basis of aerodynamic shape optimization
today~\cite{jameson1988,jameson1998}. In the geosciences, the same
construction underlies full-waveform inversion, in which the adjoint of the
wave equation converts surface seismograms into images of subsurface
structure~\cite{tarantola1984,virieux2009,fichtner2011}. In structural and
materials design, adjoint sensitivity analysis is what makes topology
optimization possible~\cite{bendsoe1988,sigmund2013}, and the same approach
now drives the inverse design of nanophotonic devices whose geometries no
human designer would propose~\cite{molesky2018}. In each case the enabling
step was not a new physical model but the ability to differentiate an
existing one efficiently.

Tutorials on the adjoint method exist across several of these communities.
Errico~\cite{Errico1997} presents the method for atmospheric and oceanic data
assimilation, Giles and Pierce~\cite{Giles2000} give a widely cited
introduction centered on aerodynamic design, and Plessix~\cite{plessix}
reviews the adjoint-state method from a geophysical perspective.
Bradley~\cite{Bradley2013} and Johnson~\cite{Johnson2012} develop
discrete-adjoint constructions in lecture-note form, and Cao, Li, Petzold, and
Serban~\cite{CaoLiPetzoldSerban2003} derive discrete adjoints for general ODE
and DAE integrators. Sirkes and Tziperman~\cite{SirkesTziperman1997} make, for
a meteorological forecast model, an observation echoed in the
discrete-versus-continuous adjoint discussion of aerodynamic
optimization~\cite{NadarajahJameson2000}: the adjoint of a discretized system
need not coincide with a discretization of the continuous adjoint equation. 

Together these treatments cover a wide range of applications and settings, continuous and discrete alike, and supply the physical context and the analytic foundations. Nevertheless, entry typically requires an advanced graduate curriculum, whether because the physics and the mathematics are developed together, or because the mathematical treatment leaves the details of an implementation for the reader to work out.

This tutorial is written for readers who have the calculus and linear algebra
needed to follow the mechanism of the adjoint method, but not yet the
physical background of a particular application or the functional analysis in
which continuous derivations are customarily presented---lower-division and
beginning graduate engineering students, and practitioners meeting the method
for the first time. It complements the discrete treatments above, which
establish the construction for a general nonlinear program or a general time
integrator, by working through a single physical inverse problem end to end,
with every vector and matrix written out and every implementation choice made
explicit. The construction reduces entirely to multivariable calculus, linear
algebra, and the transpose of a Jacobian matrix---the \emph{adjoint} of a
matrix, in the linear-algebra sense\footnote{In elementary linear algebra the
word \emph{adjoint} is sometimes used for the adjugate, or cofactor matrix.
That is a different object and plays no role here.}.

A central theme of this tutorial is that the adjoint is determined by the
discrete algorithm rather than by the differential equation itself.
Different discretizations of the same governing equation can produce adjoint
recursions that differ in their organization, update structure, and backward
propagation pattern. To make this dependence explicit, the vibrating-string
inverse problem introduced in Section~\ref{sec:inverse} is discretized in two
ways: once with displacement and velocity collocated on a single grid, and
once with velocity staggered by a half time step. For each discretization,
the forward model, its Jacobian, the adjoint equation, and the backward
recursion are derived side by side, so that the effect of the discretization
on the adjoint construction can be observed directly.

The resulting construction is closely related to reverse-mode automatic
differentiation. When the forward simulation is written as a sequence of
elementary operations, reverse-mode AD generates the same backward recursion
derived here by hand~\cite{GriewankWalther2008,Baydin2018}. The explicit
derivation nevertheless remains valuable: it exposes the underlying structure
of the gradient computation, clarifies memory and checkpointing
considerations, and provides a framework for gradient verification and
debugging.

\medskip
\noindent
\textbf{Scope.} This note omits the functional-analytic foundations, the
existence and convergence theory, and the breadth of application covered by
the standard references on optimal control, inverse problems, and
PDE-constrained
optimization~\cite{lions,plessix,gunzburger,Borzi2012,nocedal,Hinze2009,Troltzsch2010},
in favor of a single concrete construction, worked through completely.

\medskip
\noindent
\textbf{Organization.} Section~\ref{sec:inverse} introduces the inverse
problem through the vibrating string, using it to make concrete the
abstract objects---state, parameters, and governing equations---and
shows how the resulting problem takes the form of a constrained
optimization problem, $\min_p J(x,p)$ subject to $F(x,p)=0$. Section~\ref{sec:adjoint} derives the reduced gradient and the adjoint equation, and shows how the resulting gradient is verified in practice. Section~\ref{sec:backward} casts the adjoint equation as a backward recursion and specializes the construction to linear state equations. Section~\ref{sec:solver} shows how this gradient is consumed by a gradient-based
optimizer, taking L-BFGS as the working example, and examines the step
length as a requirement independent of the adjoint-computed direction.
Section~\ref{sec:discrete-wave} demonstrates the application of the
adjoint method in full detail for the vibrating string: it constructs
the two discretizations of the wave equation, and the adjoint each one
generates, once with the state variables collocated on a single grid and
once with the velocity staggered by a half time step, and closes with a
numerical example on a grid large enough to resemble practice. Section~\ref{sec:interpretation} closes the tutorial with a structured
summary of the method, a discussion of what the deterministic treatment
developed here leaves out---measurement noise and the non-uniqueness
that follows from limited data, and the Bayesian framework that
addresses both---and a closing perspective on how to think about the
adjoint method as a whole.

%%%%%%%%%%%%%%%%%%%%%%%%%%%%%%%
%%%%%%%%%%%%%%%%%%%%%%%%%%%%%%%
\section{From an Inverse Problem to Constrained Optimization}\label{sec:inverse}

Consider a string of length $L$ stretched between two fixed supports, as in
Figure~\ref{fig:string}. If the string is displaced at an interior point and
released, it vibrates in the plane determined by the two supports and the
point of release. Let $u(z,t)$ denote the transverse displacement at position
$z \in [0,L]$ and time $t$, and let $c(z)$ denote the local speed at which a
disturbance propagates along the string. Once $c(z)$ and the initial shape
are prescribed, the displacement is determined for all subsequent times:
$u(z,t)$ satisfies
\begin{equation}\label{eq:pde}
\left\{
\begin{aligned}
&\frac{\partial^{2} u}{\partial t^{2}}
   \;-\; \frac{\partial}{\partial z}
   \left( c(z)^{2}\,\frac{\partial u}{\partial z} \right) = 0,
   && z \in (0,L), \;\; t \in (0,T], \\[4pt]
&u(0,t) = 0, \qquad u(L,t) = 0,
   && t \in (0,T], \\[4pt]
&u(z,0) = u_{0}(z), \qquad
   \frac{\partial u}{\partial t}(z,0) = 0,
   && z \in [0,L].
\end{aligned}
\right.
\end{equation}
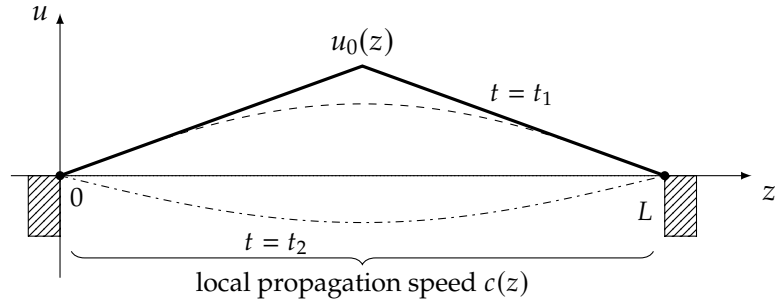
\begin{figure}[t]
\centering
\begin{tikzpicture}[>=latex]

% ---------- geometry ----------
\def\Ls{8}      % string length
\def\hp{1.45}   % pluck height

% ---------- axes ----------
\draw[->] (-0.65,0) -- (\Ls+1.15,0) node[below right] {$z$};
\draw[->] (0,-1.35) -- (0,2.15) node[left] {$u$};

% ---------- equilibrium position ----------
\draw[densely dotted] (0,0) -- (\Ls,0);

% ---------- fixed supports ----------
\fill[pattern=north east lines] (-0.42,-0.80) rectangle (0,0);
\draw (-0.42,-0.80) rectangle (0,0);
\fill[pattern=north east lines] (\Ls,-0.80) rectangle (\Ls+0.42,0);
\draw (\Ls,-0.80) rectangle (\Ls+0.42,0);

% ---------- later-time profiles ----------
\draw[dashed] plot[domain=0:\Ls,samples=120,variable=\x]
      (\x,{0.95*sin(180*\x/\Ls)});
\draw[dash pattern=on 3pt off 2pt on 1pt off 2pt]
      plot[domain=0:\Ls,samples=120,variable=\x]
      (\x,{-0.62*sin(180*\x/\Ls)});

% ---------- initial plucked shape ----------
\draw[very thick] (0,0) -- (\Ls/2,\hp) -- (\Ls,0);

% ---------- endpoints ----------
\fill (0,0) circle (1.7pt);
\fill (\Ls,0) circle (1.7pt);

% ---------- labels ----------
\node[above]        at (\Ls/2,\hp)   {$u_{0}(z)$};
\node[above right]  at (5.55,0.80)   {\small $t=t_{1}$};
\node[below right]  at (2.30,-0.62)  {\small $t=t_{2}$};

\node at (0.22,-0.28)        {\small $0$};
\node at (\Ls-0.26,-0.46)    {\small $L$};

% ---------- brace for c(z) ----------
\draw[decorate,decoration={brace,amplitude=5pt,mirror}]
      (0.15,-1.00) -- (\Ls-0.15,-1.00);
\node[below] at (\Ls/2,-1.12)
      {\small local propagation speed $c(z)$};

\end{tikzpicture}
\caption{The forward problem for a string of length $L$ held fixed at both
ends. The string is drawn aside into the shape $u_{0}(z)$ and released from
rest (solid). Once the coefficient $c(z)$ and the initial shape are
prescribed, the equation of motion \eqref{eq:pde} determines the displacement
$u(z,t)$ at every point and every later instant; two subsequent profiles are
shown.}
\label{fig:string}
\end{figure}

The first line is the wave equation in divergence form, the standard
model for the small-amplitude transverse vibration of a string whose
material properties vary along its length~\cite{Graff1991,Strauss2008}:
mass density $\rho(z)$ and tension $T(z)$ combine into the single
coefficient $c(z)=\sqrt{T(z)/\rho(z)}$, the local speed at which a
disturbance propagates, and it is the only material property that
enters the equation. This tutorial adopts~\eqref{eq:pde} as the
governing model without deriving it from a particular assumption on how
$\rho$ and $T$ individually vary; what matters for what follows is only
that $c(z)$ can vary with position, and that the resulting equation of
motion takes the divergence (self-adjoint) form
$\partial_z(c^2\partial_z u)$ rather than the naive $c^2\partial_z^2u$.
The two expressions differ whenever $c$ varies with $z$, and it is the
divergence form that must be discretized correctly.  The second line
gives the boundary conditions: the supports, shown hatched in
Figure~\ref{fig:string}, hold both ends at zero displacement for all time.
The third line gives the initial conditions: at $t=0$ the string has the
shape $u_{0}(z)$ and is at rest. For sufficiently regular $c(z)$ and
$u_{0}(z)$, system \eqref{eq:pde} has a unique solution. This is the forward
problem for the string: the material property is known and the motion is to
be determined.

The inverse problem reverses the roles of $c$ and $u$. Suppose that the
string is not uniform---thicker in some places, corroded, or unevenly
tensioned---and cannot be inspected directly. It can be excited and observed:
instruments record the motion at a limited number of locations, at a finite
number of instants, and with finite precision. The task is to recover the function $c(z)$ on all of $[0,L]$ from this
record. In the terminology of the introduction, $c(z)$ is the
\emph{parameter} of the inverse problem; because it is a function of
position, it is referred to below as the wave-speed \emph{profile}, or
simply the profile. Problems of this type
arise whenever the quantity of interest is not directly measurable, as in
seismology, where ground motion recorded at the surface is used to infer the
structure of the earth below. In each case the unknown is a property
distributed through a medium, the observable is a partial sampling of a
field governed by that property, and the two are related by a partial
differential equation such as \eqref{eq:pde}. The unknown does not appear in
the data explicitly: the recorded motion depends on $c(z)$ only through the
solution of \eqref{eq:pde}. A candidate profile can therefore be assessed only by simulating the string
and comparing the predicted record with the measured one. Three properties of
the problem determine how it can be solved. The problem is nonlinear: the
equation of motion is linear in the displacement, but the map from $c$ to the
recorded data is not, because $c$ enters as a coefficient of the differential
operator rather than as a forcing term, and scaling $c$ by a constant does not
scale the recorded data by that constant. The problem is ill-posed: the
measurements are finite in number and contaminated by noise, whereas the
unknown is a function of position, so many distinct profiles reproduce the
record to within the noise level, small perturbations of the record may
correspond to large changes in the inferred profile, and no unique, stable
inverse exists. The problem is large: once discretized, the unknown profile
may have thousands to millions of degrees of freedom, and any method whose
cost grows in proportion to that number is impractical.

The first two difficulties are addressed by replacing the search for an
exact inverse with the search for the profile that best accounts for the
measurements. This requires a scalar measure of the failure of a candidate
profile. Let $z_{1},\dots,z_{n_{s}}$ denote the sensor locations,
$t_{1},\dots,t_{n_{r}}$ the recording instants, and $d_{s,j}$ the
displacement recorded at sensor $s$ at instant $t_{j}$. The discrepancy
between prediction and observation is measured by
\begin{equation}\label{eq:misfit}
J \;=\; \frac{1}{2}\sum\nolimits_{j=1}^{n_{r}}\sum\nolimits_{s=1}^{n_{s}}
\big( u(z_{s},t_{j}) - d_{s,j} \big)^{2}
\;+\; \frac{\beta}{2}\,\mathcal{R}(c),
\end{equation}
whose minimizer is sought. The outer sum runs over recording instants and the
inner sum over sensors; $u(z_{s},t_{j})$ is the solution of \eqref{eq:pde}
evaluated at the sensor location and recording instant, so $J$ depends on
$c$ only through these sampled values. The first term penalizes disagreement
with the measurements. The second is a regularization term, weighted by
$\beta > 0$, that restores well-posedness by expressing a preference among
the profiles that fit the data equally well; a penalty on roughness, for
example, favors profiles that vary smoothly along the string. The particular
form of \eqref{eq:misfit} plays no role in what follows; only
differentiability is required.

This reformulation introduces the distinction on which the entire
development rests. Two objects appear in the objective: the profile, which
may be adjusted freely, and the motion of the string, which may not. If the
displacement were treated as an independent unknown, the mismatch could be
driven to zero by a displacement field that is not a solution of
\eqref{eq:pde}; the governing system must therefore be imposed as a
\emph{constraint} on the minimization, not merely used to generate
predictions. Writing \eqref{eq:pde} in residual form---so that satisfaction
of the equation, together with its boundary and initial conditions, is
equivalent to the vanishing of a residual---casts the inverse problem as:
minimize the mismatch over all admissible parameters, subject to the
requirement that the state satisfy the governing equations. This is a
\emph{PDE-constrained optimization problem}.

This problem cannot be solved directly: since $J$ depends on the parameters
only through the solution of a constraint it does not appear in explicitly,
its stationarity conditions form a nonlinear system no easier to solve than
the inverse problem itself. It is instead solved iteratively. Given a value
of the parameters, the forward problem is solved, the objective is
evaluated, and a new value is chosen at which the objective is smaller;
repeating this step produces a sequence of parameter values, the
\emph{iterates}, along which the objective decreases monotonically and,
being bounded below, converges---under standard assumptions, to a
stationary point. A method of this type is a \emph{successive descent
method}: rather than solving for the minimizer directly, it approaches the
minimizer through a sequence of improvements, each requiring one solution of
the forward problem, with the corresponding state satisfying the constraint
exactly at every step and the iteration itself taking place in the
parameters~alone.

Each step of a descent method consists of a \emph{descent direction}---a
direction in parameter space along which the objective decreases for
sufficiently small steps---and a \emph{step size} specifying how far to move
along it. The step size is determined by line-search procedures that are
standard and not specific to this problem. The direction is determined by the
first-order behavior of the objective. For a differentiable objective, the
change produced by a small displacement of the parameters is, to first order,
the inner product of the displacement with the gradient of the objective with
respect to the parameters. A direction is a descent direction if and only if
this inner product is negative, and the negative gradient is the direction of
steepest decrease. The gradient therefore determines which directions are
descent directions, and computing it is the central difficulty of the
problem: the dependence noted above runs only through the state, so
differentiating the objective means differentiating through the governing
equations. The adjoint method is the technique that carries out this
differentiation efficiently.

One decision precedes its construction. The governing equations are
continuous, but all computation is performed on a discrete system. One may
differentiate first and discretize the resulting adjoint equations, or
discretize first and differentiate the finite-dimensional system that
results. These routes are called \emph{optimize-then-discretize} and
\emph{discretize-then-optimize}, respectively. They do not in general yield
the same numerical gradient, because differentiation and discretization do
not commute. This tutorial follows the second route, for the reasons given in
the introduction: the derivation requires only calculus and linear algebra,
and the gradient it produces is exactly the gradient of the discretized
objective that a program evaluates.

Approximating \eqref{eq:pde} on a grid in space and time reduces it to a
finite system of algebraic equations, and the constrained problem studied in
the remainder of this tutorial is
\begin{equation}\label{eq:opt}
\min_{p} J(x, p) \quad \text{subject to} \quad F(x, p) = 0,
\end{equation}
where
\begin{itemize}[leftmargin=2em]
  \item $x \in \mathbb{R}^{n}$ is the vector of \textbf{state variables},
  \item $p \in \mathbb{R}^{m}$ is the vector of \textbf{parameters},
  \item $F : \mathbb{R}^{n} \times \mathbb{R}^{m} \to \mathbb{R}^{n}$ is the
  discrete \textbf{state equation}.
\end{itemize}
Equation~\eqref{eq:opt} is not specific to the vibrating string: it is the
form taken by any PDE-constrained inverse problem once the governing
equations have been discretized, with $F$ encoding whichever equation and
discretization scheme are in use. The string is simply the concrete instance
carried through this tutorial, chosen because it is the simplest setting in
which every object above can be made explicit. Each of these objects, and the
objective $J$, is now identified for that instance. The abstract notation of
\eqref{eq:opt} is the point at which first-time readers most often lose track
of what the symbols denote.

\medskip
\noindent
\textbf{The parameters $p$.} The unknown in the inverse problem is the
continuous profile $c(z)$. Discretizing the spatial domain into $n_{z}$
interior nodes replaces this function by the finite collection of its values
on the resulting grid, and this collection is the parameter vector
$p \in \mathbb{R}^{m}$. Two choices are left open here and settled by the
discretization. The first is whether to take $c$ or $c^{2}$ as the parameter;
since $c$ enters \eqref{eq:pde} only through $c^{2}$, the squared values are
often more convenient. The second is where the parameter values are located:
a scheme that evaluates the coefficient at the nodes gives $m = n_{z}$,
whereas a scheme that evaluates it between nodes gives $m = n_{z}+1$.
Sections~\ref{sec:collocated} and~\ref{sec:staggered} make one choice each,
and the gradients they produce are therefore gradients with respect to
different parameter vectors.

\medskip
\noindent
\textbf{The state $x$.} The wave equation is second order in time, so the
state at a single instant is not the displacement alone; the subsequent
motion is determined by the displacement together with its rate of change.
Setting $v = \partial u/\partial t$ and sampling both fields at the $n_{z}$
spatial nodes, the state at time level $k$ is
\[
x^{k} =
\begin{bmatrix} u^{k} \\ v^{k} \end{bmatrix}
\in \mathbb{R}^{2 n_{z}},
\qquad
u^{k} = \big[u^{k}_{1}, \dots, u^{k}_{n_{z}}\big]^{\top},
\quad
v^{k} = \big[v^{k}_{1}, \dots, v^{k}_{n_{z}}\big]^{\top},
\]
and the full state vector stacks these over all $n_{t}$ time levels,
\[
x = \big[\, (x^{1})^{\top},\; (x^{2})^{\top},\; \dots,\;
(x^{n_{t}})^{\top} \,\big]^{\top} \in \mathbb{R}^{n},
\qquad n = 2\, n_{z}\, n_{t}.
\]
The definition of the state depends on the discretization. A scheme that
places $u$ and $v$ on the same time levels and one that offsets them by a
half step lead to different, equally valid, definitions of $x^{k}$.
Sections~\ref{sec:collocated} and~\ref{sec:staggered} carry out both
constructions and show that the choice affects the adjoint.

\medskip
\noindent
\textbf{The state equation $F$.} The discretized equation of motion advances
the state from one time level to the next. The boundary conditions at the two
supports are built into the discrete spatial operator, and the prescribed
initial state $x^{0}$ enters the first update as known data, not as an
unknown; it is not a component of $x$, and the consequence of this for the
adjoint is noted in Section~\ref{sec:collocated}. Collecting the update rule
at every time level and writing each as a residual required to vanish gives a
single system of $n$ algebraic equations,
\[
F(x,p) = 0,
\]
whose solution is the entire simulated history. Running the simulation and
solving $F(x,p)=0$ are the same operation.

\medskip
\noindent
\textbf{The objective $J$.} After discretization the sensors occupy a set
$\mathcal{S}\subset\{1,\dots,n_{z}\}$ of grid nodes, with
$|\mathcal{S}| = n_{s}$, and the recording instants are taken to coincide
with the time levels, so that $n_{r} = n_{t}$. The sampled solution
$u(z_{s},t_{j})$ in \eqref{eq:misfit} is then the state component
$u^{k}_{i}$ with $i \in \mathcal{S}$ and $k = j$, the recorded value
$d_{s,j}$ is written $d^{k}_{i}$, and the misfit becomes
\begin{equation}\label{eq:objective}
J(x,p) \;=\; \frac{1}{2}\sum\nolimits_{k=1}^{n_{t}}\sum\nolimits_{i\in\mathcal{S}}
\big( u^{k}_{i} - d^{k}_{i} \big)^{2}
\;+\; \frac{\beta}{2}\,\mathcal{R}(p).
\end{equation}
Let $d$ denote the vector of all $d^{k}_{i}$ in a fixed ordering, and let
$H$ be the matrix that extracts from $x$ the components $u^{k}_{i}$,
$i\in\mathcal{S}$, in the same ordering. Then $Hx$ is the vector of
predicted data, the double sum equals $\|Hx-d\|^{2}$, and
\[
J(x,p) \;=\; \tfrac{1}{2}\,\big\|Hx-d\big\|^{2}
\;+\; \tfrac{\beta}{2}\,\mathcal{R}(p),
\qquad
\frac{\partial J}{\partial x} \;=\; H^{\top}(Hx-d),
\qquad
\frac{\partial J}{\partial p} \;=\; \tfrac{\beta}{2}\,\nabla\mathcal{R}(p).
\]
The compact form is used in the derivations that follow; the double sum is
what a program evaluates. The data term depends on $x$ and not on $p$, and
$\partial J/\partial x$ is nonzero only in the observed components; the
regularization term depends on $p$ and not on $x$. The coupling between
$x$ and $p$ enters only through the constraint $F(x,p)=0$. If sensors do not
lie on grid nodes or recording instants do not coincide with time levels,
$H$ becomes an interpolation matrix rather than a selection matrix, and
nothing in the development changes.

Before proceeding further, it is worth mentioning that nothing in the
constructions of Sections~\ref{sec:adjoint}, \ref{sec:backward},
and~\ref{sec:solver} depends on the particular form of \eqref{eq:pde} or
\eqref{eq:opt}; they are stated for the general state equation $F(x,p)=0$ and
apply unchanged to any discretized PDE-constrained inverse problem, the
vibrating string being only the vehicle used to make each step concrete.
%%%%%%%%%%%%%%%%%%%%%%%%%%%
%%%%%%%%%%%%%%%%%%%%%%%%%%%
\section{Adjoint Method}\label{sec:adjoint}
The inverse problem developed in the previous section leads to the
PDE-constrained optimization problem \eqref{eq:opt}, solved by the
successive descent methods of Section~\ref{sec:inverse}: at every
iterate, a descent direction is required, and computing it is the
central difficulty, since $J$ depends on the parameters only through the
state. The purpose of this section is to construct that direction
explicitly, by examining how the constraint $F(x,p)=0$ shapes it.

Viewed purely as a constrained optimization problem, \eqref{eq:opt}
contains two sets of decision variables, the state $x$ and the parameter
$p$, and a conventional treatment would consider perturbing both when
constructing a descent direction. But in a PDE-constrained optimization these two variables do not play the same
role: as established above, $p$ is free while $x$ is determined by it
through $F(x,p)=0$, so a perturbation of $p$ forces a corresponding
perturbation of $x$ rather than leaving it free to choose independently.
Exploiting this dependence is the key observation underlying the
adjoint method.

To arrive at that construction, we take a step back and revisit the basic mechanism of iterative descent methods. The discussion that follows may appear to be a detour into elementary optimization, but it is precisely here that the adjoint method originates. A clear understanding of these steps removes much of the mystery surrounding the adjoint derivation and makes its role in PDE-constrained optimization transparent.

Starting from an initial parameter estimate $p_0$, an iterative method generates a sequence 
$$p_{\ell+1}=p_\ell+\alpha_\ell d_\ell,\quad \ell\in\{0,1,2,\cdots\}$$
where $d_\ell$ is a descent direction and $\alpha_\ell$ is a step size. The fundamental question is therefore: how should $d_\ell$ and $\alpha_\ell$ be chosen so that the objective decreases while the governing equations remain satisfied? To answer this question, we examine how the objective changes under small perturbations of the state and parameter variables.

Let $(x_\ell,p_\ell)$ be a feasible point, i.e., \[ F(x_\ell,p_\ell)=0. \] To determine whether a modification of the parameter reduces the objective, we consider a small perturbation $\Delta p$ of the current parameter value. Since the state and parameter are linked through the constraint equation, this perturbation induces a corresponding state perturbation $\Delta x$. The resulting change in the objective can therefore be understood by examining the simultaneous perturbation \[ (x_\ell,p_\ell)\;\mapsto\; (x_\ell+\Delta x,\;p_\ell+\Delta p). \] 
The effect of these perturbations on the objective is obtained from the first-order Taylor expansion \begin{equation} J(x_\ell + \Delta x, p_\ell + \Delta p) = J(x_\ell, p_\ell) + \left( \frac{\partial J}{\partial x} \right)^\top \Delta x + \left( \frac{\partial J}{\partial p} \right)^\top \Delta p + \mathcal{O}(\|\Delta x,\Delta p\|^2). \end{equation}
The gradients are defined as column vectors:
\begin{equation}
\frac{\partial J}{\partial x} \in \mathbb{R}^{n \times 1} = \begin{bmatrix} \frac{\partial J}{\partial x_1} \\ \vdots \\ \frac{\partial J}{\partial x_n} \end{bmatrix}, \qquad \frac{\partial J}{\partial p} \in \mathbb{R}^{m \times 1} = \begin{bmatrix} \frac{\partial J}{\partial p_1} \\ \vdots \\ \frac{\partial J}{\partial p_m} \end{bmatrix}.
\end{equation}

Neglecting higher-order terms, a perturbation decreases the objective,
i.e., $$J(x_{\ell}+\Delta x,\,p_{\ell}+\Delta p) \le J(x_{\ell},p_{\ell}),$$
whenever the first-order variation is non-positive. This gives the
\emph{descent condition}
\begin{equation}\label{eq:descent}
\left(\frac{\partial J}{\partial x}\right)^{\top} \Delta x
\;+\;
\left(\frac{\partial J}{\partial p}\right)^{\top} \Delta p
\;\le\; 0.
\end{equation}

The descent condition \eqref{eq:descent} involves both $\Delta x$ and
$\Delta p$, but, as already noted, these variations are not independent:
once $\Delta p$ is chosen, $\Delta x$ is fixed as the state response induced
by it through the constraint equation. Making this dependence explicit is
therefore the next step. Since both $(x_\ell,p_\ell)$ and
$(x_\ell+\Delta x,\,p_\ell+\Delta p)$ must satisfy the governing equations,
\begin{equation}
F(x_\ell+\Delta x,\,p_\ell+\Delta p)=0.
\end{equation}
Applying a first-order Taylor expansion gives
\begin{equation}
F(x_\ell, p_\ell)
+
\left( \frac{\partial F}{\partial x} \right)_{n \times n}
(\Delta x)_{n \times 1}
+
\left( \frac{\partial F}{\partial p} \right)_{n \times m}
(\Delta p)_{m \times 1}
+ \text{H.O.T.}
= 0 .
\end{equation}
where \[\frac{\partial F}{\partial x}=\begin{bmatrix}\frac{\partial F_1}{\partial x_1}&\frac{\partial F_1}{\partial x_2}\cdots&\frac{\partial F_1}{\partial x_n}\\    
\frac{\partial F_2}{\partial x_1}&\frac{\partial F_2}{\partial x_2}\cdots&\frac{\partial F_2}{\partial x_n}\\    
\vdots\\
\frac{\partial F_n}{\partial x_1}&\frac{\partial F_n}{\partial x_2}\cdots&\frac{\partial F_n}{\partial x_n}\\    
\end{bmatrix},\qquad \frac{\partial F}{\partial p}=\begin{bmatrix}\frac{\partial F_1}{\partial p_1}&\frac{\partial F_1}{\partial p_2}\cdots&\frac{\partial F_1}{\partial p_m}\\    
\frac{\partial F_2}{\partial p_1}&\frac{\partial F_2}{\partial p_2}\cdots&\frac{\partial F_2}{\partial p_m}\\    
\vdots\\
\frac{\partial F_n}{\partial p_1}&\frac{\partial F_n}{\partial p_2}\cdots&\frac{\partial F_n}{\partial p_m}\\    
\end{bmatrix}
\]
Since $F(x_\ell,p_\ell)=0$, the variation dependency is
\begin{equation}
\frac{\partial F}{\partial x} \Delta x
+
\frac{\partial F}{\partial p} \Delta p
= 0
\quad \Longrightarrow \quad
\Delta x
=
-
\left( \frac{\partial F}{\partial x} \right)^{-1}
\frac{\partial F}{\partial p}
\Delta p .
\end{equation}
We assume that $\frac{\partial F}{\partial x}$ is invertible, which corresponds to the forward problem being well-posed. Substituting this relation into the descent condition eliminates $\Delta x$ and expresses the first-order variation of the objective entirely in terms of the parameter perturbation:
\begin{equation}
\left( -\left(\frac{\partial F}{\partial p}\right)^\top
\left(\frac{\partial F}{\partial x} \right)^{-\top}
\frac{\partial J}{\partial x}
+
\frac{\partial J}{\partial p}
\right)^\top
\Delta p
\le 0 .
\end{equation}

Taking $\Delta p$ as free variation, a sufficient choice for obtaining first-order descent is therefore 
\begin{equation}\label{eq::descent}
d_\ell=\Delta p = -\left( -\left(\frac{\partial F}{\partial p}\right)^{\top} \left(\frac{\partial F}{\partial x}\right)^{-\top} \frac{\partial J}{\partial x} + \frac{\partial J}{\partial p} \right). \end{equation} 
This expression determines how the parameters should be modified to reduce the objective while respecting the constraint.
 
It is instructive to compare this expression with the derivative of the \emph{reduced objective} \[ \hat J(p)=J(x(p),p), \] obtained by regarding the state as an implicit function of the parameters through the constraint equation \(F(x(p),p)=0\). By the chain rule, \[ \frac{d\hat J}{dp} = \left(\frac{\partial x}{\partial p}\right)^\top \frac{\partial J}{\partial x} + \frac{\partial J}{\partial p}. \]

Differentiating the constraint \[ F(x(p),p)=0 \] with respect to \(p\) gives \[ \frac{\partial F}{\partial x} \frac{\partial x}{\partial p} + \frac{\partial F}{\partial p} = 0, \] and therefore \[ \frac{\partial x}{\partial p} = -\left(\frac{\partial F}{\partial x}\right)^{-1} \frac{\partial F}{\partial p}. \]

Substituting this relation into the chain rule yields 
\begin{equation}\label{eq:reduced_grad1}
\frac{\text{d}\hat J}{\text{d}p} = -\left(\frac{\partial F}{\partial p}\right)^{\top} \left(\frac{\partial F}{\partial x}\right)^{-\top} \frac{\partial J}{\partial x} + \frac{\partial J}{\partial p}, 
\end{equation} which is precisely the quantity appearing in the descent direction \eqref{eq::descent}. We therefore identify 
\begin{equation}\label{eq::g_reduced}\nabla_p \hat J = -\left(\frac{\partial F}{\partial p}\right)^{\top} \left(\frac{\partial F}{\partial x}\right)^{-\top} \frac{\partial J}{\partial x} + \frac{\partial J}{\partial p} 
\end{equation}
as the gradient of the reduced objective.
Since the optimization is carried out entirely in the parameter space, we henceforth suppress the hat and denote the reduced objective simply by \(J(p)\). The descent direction is therefore \[ d_\ell = -\nabla_p J(p_\ell). \]

Thus, the descent direction is the negative gradient of the reduced objective. At first sight, this appears to solve the optimization problem: the gradient has been obtained and a corresponding descent direction follows immediately. The resulting expression is, however, not yet a practical computational procedure. Interpreted literally, \eqref{eq::g_reduced} involves the inverse transpose Jacobian \( \left(\frac{\partial F}{\partial x}\right)^{-\top}, \) where \(\partial F/\partial x\) is the Jacobian of the discretized governing equations. For large-scale PDE discretizations, this matrix has dimension equal to the number of state $x$ and may be extremely large. The key computational question is therefore not how to form the inverse, but how to evaluate the reduced-gradient formula efficiently.

%%%%%

One possibility is to recognize that the factor $\left(\frac{\partial F}{\partial x}\right)^{-1}\frac{\partial F}{\partial p}$ appearing in \eqref{eq::g_reduced} is exactly the state sensitivity matrix $\frac{\partial x}{\partial p}\in\mathbb{R}^{n\times m}$, obtained by differentiating the constraint $F(x,p)=0$ directly with respect to $p$:\footnote{Differentiating $F(x,p)=0$ with respect to $p$ gives $\frac{\partial F}{\partial x}\frac{\partial x}{\partial p}+\frac{\partial F}{\partial p}=0$, so that $\frac{\partial x}{\partial p}=-\left(\frac{\partial F}{\partial x}\right)^{-1}\frac{\partial F}{\partial p}$.}
\begin{equation*}
\nabla_p \hat J =\underbrace{ -\left(\frac{\partial F}{\partial p}\right)^{\top} \left(\frac{\partial F}{\partial x}\right)^{-\top}}_{\left(\frac{\partial x}{\partial p}\right)^\top} \frac{\partial J}{\partial x} + \frac{\partial J}{\partial p} .
\end{equation*}
Column $j$ of this matrix, $s_j=\partial x/\partial p_j\in\mathbb{R}^{n}$, is itself a vector of the same dimension as the state, obtained by solving
\begin{equation}
\frac{\partial F}{\partial x}\,s_j=-\frac{\partial F}{\partial p_j},\qquad j=1,\dots,m.
\end{equation}
Computing $\partial x/\partial p$ this way, which is known as the \emph{forward-sensitivity method}, therefore requires $m$ separate $n$-dimensional linear solves, one for each parameter, each with the same coefficient matrix $\partial F/\partial x$ but a different right-hand side.

The more computationally efficient path, however, is to read the first summand of \eqref{eq::g_reduced} from right to left: first form the vector $\partial J/\partial x$, then solve a single linear system with the transpose coefficient matrix, and only afterward contract the result against $\partial F/\partial p$. Concretely, introduce the vector $\lambda\in\mathbb{R}^n$, defined as the solution of
\begin{equation}\label{eq::adjoint_eq}
\underbrace{-\left(\frac{\partial F}{\partial x}\right)^{-\top}\frac{\partial J}{\partial x}}_{\lambda} \quad\Longleftrightarrow\quad \left(\frac{\partial F}{\partial x}\right)^{\top}\lambda=-\frac{\partial J}{\partial x},
\end{equation}
so that the reduced gradient becomes
\begin{equation}
\nabla_p J=\left(\frac{\partial F}{\partial p}\right)^{\top}\lambda+\frac{\partial J}{\partial p}.
\label{eq:reduced_grad2}
\end{equation}
The matrix $\left(\dfrac{\partial F}{\partial x}\right)^{\top}$ is the
transpose, or adjoint, of the Jacobian $\dfrac{\partial F}{\partial
x}$. For this reason, the system~\eqref{eq::adjoint_eq} above is called the \emph{adjoint
equation}, the variable $\lambda$ is called the \emph{adjoint
variable}, and the resulting gradient-computation procedure is known as
the \emph{adjoint method}. 

Both the forward-sensitivity method and the adjoint method are exact, and neither ever
inverts $\partial F/\partial x$ as a matrix; what distinguishes them is
solely the number of linear solves, $m$ against $1$---the same
associativity that separates forward-mode from reverse-mode automatic
differentiation---and it is decisive whenever the number of parameters
$m$ is large, which for a discretized field is the rule rather than the
exception.

The vector $\partial J/\partial x$ measures how perturbations of the
state variables affect the objective function. The role of the adjoint
variable $\lambda$ is to carry this sensitivity information from the
state space to the parameter space without ever forming the state
sensitivity matrix $\partial x/\partial p$ that the forward-sensitivity
route would require.

%Substituting the adjoint variable into the reduced-gradient expression gives a computational procedure that requires only the solution of the adjoint system and the evaluation of the remaining derivatives. 

%%%%%%%%%%%%%%%%%%%%%%%%%%%
\section{The Adjoint as a Backward Recursion}
\label{sec:backward}

The adjoint equation~\eqref{eq::adjoint_eq}, at first encounter, reads as a single, large-scale linear system. One of the advantages of the discrete formulation, however, is that it exposes the structural dependency among the time-marched variables: for time-dependent problems, solving this system reduces to a backward march, and the adjoint equation is naturally solved backward in time in an iterative manner.

To see why, consider a generic time-marching scheme in which the state
at time level $k+1$ is obtained from the state at time level $k$
through
\[
F_k(x_{k+1},x_k,p)=0, \qquad k=0,\ldots,N-1,
\]
with the initial state $x_0$ prescribed as data, exactly as in
Section~\ref{sec:inverse}.

The unknown state history is the stacked vector
\[
x = \begin{bmatrix} x_1 \\ x_2 \\ \vdots \\ x_N \end{bmatrix},
\]
and the $N$ governing equations $F_0,\dots,F_{N-1}$ assemble into the
single residual system $F(x,p)=0$: an $N$-block system for the $N$
unknowns $x_1,\dots,x_N$. The initial state $x_0$ enters as data in the
first block equation, $F_0(x_1,x_0,p)=0$, and is not itself a
component of $x$.

Because $F_k$ depends only on $x_{k+1}$ and $x_k$, the Jacobian
$\partial F/\partial x$ with respect to the unknowns $x_1,\dots,x_N$
has the block lower-triangular structure
\[
\frac{\partial F}{\partial x} =
\begin{bmatrix}
B_0 & 0 & 0 & \cdots & 0 \\
A_1 & B_1 & 0 & \cdots & 0 \\
0 & A_2 & B_2 & \cdots & 0 \\
\vdots & & \ddots & \ddots & \vdots \\
0 & \cdots & 0 & A_{N-1} & B_{N-1}
\end{bmatrix},
\]
an $N\times N$ block matrix with $B_k=\partial F_k/\partial x_{k+1}$ and
$A_k=\partial F_k/\partial x_k$ for $k=1,\dots,N-1$; the first row has
only $B_0=\partial F_0/\partial x_1$, since the dependence of $F_0$ on
the data $x_0$ contributes no column---$x_0$ is not among the
unknowns. This lower-triangular structure expresses the causal nature
of the forward simulation: information propagates from earlier times
to later times.

The adjoint equation involves the transpose of this Jacobian,
\[
\left(\frac{\partial F}{\partial x}\right)^\top =
\begin{bmatrix}
B_0^\top & A_1^\top & 0 & \cdots & 0 \\
0 & B_1^\top & A_2^\top & \cdots & 0 \\
\vdots & & \ddots & \ddots & A_{N-1}^\top \\
0 & \cdots & 0 & 0 & B_{N-1}^\top
\end{bmatrix}.
\]

The transpose converts the lower-triangular structure of the forward
problem into an upper-triangular structure. Consequently, the adjoint
variables cannot be computed in chronological order. Instead, they
must be determined starting from the final time level and proceeding
backward.

To make this explicit, write the adjoint vector as
\[
\lambda = \begin{bmatrix} \lambda_1 \\ \lambda_2 \\ \vdots \\ \lambda_N \end{bmatrix},
\]
and similarly partition
\[
-\frac{\partial J}{\partial x} = \begin{bmatrix} r_1 \\ r_2 \\ \vdots \\ r_N \end{bmatrix},
\]
both with respect to the unknowns $x_1,\dots,x_N$ alone; since $x_0$ is
data rather than an unknown, no corresponding block $\lambda_0$ or
$r_0$ exists. The adjoint equation
\[
\left(\frac{\partial F}{\partial x}\right)^\top \lambda = -\frac{\partial J}{\partial x}
\]
then expands into
\[
\begin{aligned}
B_{N-1}^\top \lambda_N &= r_N,\\
B_{N-2}^\top \lambda_{N-1} + A_{N-1}^\top \lambda_N &= r_{N-1},\\
B_{N-3}^\top \lambda_{N-2} + A_{N-2}^\top \lambda_{N-1} &= r_{N-2},\\
&\vdots
\end{aligned}
\]
The final adjoint variable $\lambda_N$ is obtained first. Once
$\lambda_N$ is known, the equation above determines $\lambda_{N-1}$,
then $\lambda_{N-2}$, and so on until $\lambda_1$ is reached; there is
no $\lambda_0$ to solve for, since $x_0$ is prescribed data rather than
an unknown of the forward problem. The adjoint therefore satisfies a
backward recursion.

\paragraph{The terminal data $r_N$ and the role of boundary conditions}
The value of $r_N$ is fixed, once and for all, by the objective:
$r_N = -\partial J/\partial x\big|_{x^N}$, the negative gradient of $J$
with respect to whichever components of the final state $x^N$ it
depends on. Nothing about the history of the state or the choice of
boundary condition enters this formula; it is determined only by which
entries of $x^N$ the objective observes. For the quadratic misfit~\eqref{eq:objective}, $J(x,p) =
\tfrac12\|Hx-d\|^2 + \tfrac{\beta}{2}\mathcal{R}(p)$, so
$\partial J/\partial x = H^\top(Hx-d)$ and
$r_N = -H_N^\top(H_N x^N - d_N)$, a vector that is zero in every
component of $x^N$ not sampled by a sensor at the final recording
instant. The regularization term contributes nothing to $r_N$: since
$\mathcal{R}(p)$ depends on $p$ alone, it enters the reduced gradient
only through $\partial J/\partial p = \tfrac{\beta}{2}\nabla\mathcal{R}(p)$
in $\nabla_p J = \partial J/\partial p + (\partial F/\partial p)^\top
\lambda$, never through the state-side terminal data that drives the
backward recursion.

What the boundary conditions determine is not this formula but the
meaning of $x^N$ itself---which physical degrees of freedom are
components of the state at all. A boundary condition that prescribes a
quantity outright, rather than solving for it, removes the corresponding
degree of freedom from $x^k$ at every time level: it is data supplied to
$F(x,p)=0$, not an unknown determined by it, so it never appears in
$x^N$, and consequently no component of $r_N$, and no adjoint variable,
is ever associated with it---the same mechanism already noted above for
the prescribed initial state $x^0$. This is the case throughout
Sections~\ref{sec:discrete-wave}, where both ends of the string are held
fixed, $u(0,t) = u(L,t) = 0$, and the boundary displacements are exactly
such prescribed data.

Suppose instead that the end $z=L$ is left free while a disturbance is
launched from the fixed end at $z=0$. A free end imposes the Neumann
condition $\partial u/\partial z(L,t)=0$, and a discretization of it
typically replaces the elimination used at a fixed end with a one-sided
or ghost-node stencil, leaving the boundary displacement $u_{n_z}^k$ a
genuine unknown, present in $x^k$ for every $k$. Two consequences
follow, and only one of them touches $r_N$. First, $u_{n_z}^N$ is now a
component of $x^N$, so $r_N$ acquires a corresponding entry---nonzero
only if a sensor happens to observe that node at the final instant, as
for any newly admitted degree of freedom; the free boundary condition
itself contributes nothing to $r_N$ beyond making the slot exist.
Second, and this is where the boundary condition actually acts on the
adjoint, the blocks $A_k$ and $B_k$ built at that end use the free-end
stencil rather than the fixed-end elimination, so $A_k^\top$ and
$B_k^\top$ carry different coefficients there than they would for a
fixed end. Since every backward step multiplies by these transposed
blocks, this---not a change to the terminal condition---is what governs
how sensitivity launched at $z=0$ is reflected once it reaches the free
end at $z=L$.

This is the discrete counterpart of a fact familiar from the calculus of
variations, where a Dirichlet condition forces the corresponding
variation to vanish at the boundary while a Neumann condition instead
leaves a boundary term, produced by integration by parts, that must
vanish on its own---a natural boundary condition on the adjoint, derived
separately for each type of primal condition. No such separate
derivation is needed here: the discretization of the boundary condition
is built into $A(p)$ once, for whatever end conditions are chosen, and
its transpose, used verbatim in the recursion above, already encodes
whatever the correct adjoint behavior at that boundary is. Changing the
boundary condition changes the forward operator; it does not change the
rule by which the adjoint is derived from it.

The reversal from a lower-triangular forward system to an
upper-triangular adjoint system explains several features that often
appear mysterious in continuous adjoint derivations:
\begin{itemize}
\item The forward state is propagated forward because the governing
equations are lower triangular in time.
\item The adjoint is propagated backward because the transpose of a
lower triangular matrix is upper triangular.
\item The terminal condition of the adjoint appears naturally because
the backward recursion begins at the final time level.
\item No adjoint variable is associated with prescribed initial
conditions, since those quantities are not unknowns of the forward
problem.
\end{itemize}

The key point is that backward propagation is not an additional physical principle. It is a consequence of linear algebra. Once the governing equations have been discretized in time, the transpose appearing in the adjoint equation reverses the direction in which information must be propagated through the computational graph. This is the origin of the backward-in-time adjoint recursions encountered throughout PDE-constrained optimization.

\paragraph{Storage and the cost of the backward sweep}
The reduced gradient's cost is often summarized as ``one forward solve
plus one adjoint solve,'' and in the number of \emph{state-equation
solves} this is exactly right: computing $\nabla_p J$ requires forming
$x^0,\dots,x^N$ once, forward, and $\lambda^N,\dots,\lambda^1$ once,
backward---never one solve per parameter. What this count leaves out is
memory. The backward recursion needs $\partial F/\partial x$ evaluated
at $x^k$, and the parameter-gradient assembly needs $x^k$ itself, at
every level the backward sweep visits; since that sweep runs from
$k=N$ down to $k=0$ only after the forward sweep has already reached
$k=N$, every intermediate state must either be stored or recomputed
when the backward sweep later needs it. Storing the full trajectory
costs $O(n\,n_t)$ memory, with $n$ the state dimension and $n_t$ the
number of time levels---negligible for the small examples worked out
by hand in this tutorial, and even for the numerical demonstration of
Section~\ref{sec:numerical-demonstration}, but often the binding
constraint in practice: a three-dimensional simulation with $n$ in the
millions and $n_t$ in the thousands can make storing every time level
infeasible well before the forward solve itself becomes prohibitive.
The standard remedy is \emph{checkpointing}: store the state at a
sparse set of time levels during the forward sweep, and, whenever the
backward sweep needs a state that was not stored, recompute it by
re-running the forward recursion from the nearest earlier checkpoint.
Griewank and Walther's \texttt{revolve}
algorithm~\cite{GriewankWalther2000,GriewankWalther2008} allocates a
fixed checkpoint budget optimally, trading $O(\log n_t)$ storage for
$O(n_t\log n_t)$ forward re-solves instead of $O(n_t)$ storage at
$O(n_t)$ solves, and some form of it is implemented in essentially
every adjoint or automatic-differentiation framework used at scale.
The ``one forward, one adjoint'' cost statement is therefore a
statement about solves, not about memory, and the two need not be
minimized together.

\subsection{Linear State Equations}
The general adjoint construction of Section~\ref{sec:adjoint} required
differentiating $F$ with respect to $x$ and $p$ in the abstract, an
operation that has to be redone for every new governing equation. The
case treated in this subsection removes that burden for a large and
practically important class of problems: whenever the discretized state
equation is linear in $x$, the adjoint operator is available for free,
as the transpose of a matrix the forward solver has already assembled,
with no separate differentiation step at all. This is not a special
case chosen for mathematical convenience---it is the case that arises
whenever a PDE is discretized by a linear time-marching scheme, which
includes both discretizations of the wave equation developed in
Section~\ref{sec:discrete-wave}. Working through it here in general
form means that section can proceed directly to writing down the
adjoint recursion, without re-deriving it from scratch.

Many discretized PDEs lead to a state equation that is linear in the
state variables. In this case the residual takes the form
\[
F(x,p)=A(p)x+b(p),
\]
where the matrix \(A(p)\) and the vector \(b(p)\) may depend on the
parameters \(p\). The adjoint equation then becomes particularly
transparent.

To construct it, we require the Jacobians \(\partial F/\partial x\) and
\(\partial F/\partial p\).

\paragraph{Derivative with respect to the state}
Since \(A(p)\) and \(b(p)\) depend only on the parameters and not on
the state, differentiation with respect to \(x\) gives
\[
\frac{\partial F}{\partial x} = A(p).
\]
Consequently, the adjoint equation
\[
\left(\frac{\partial F}{\partial x}\right)^T\lambda = -\frac{\partial J}{\partial x}
\]
reduces to
\[
A(p)^T\lambda = -\frac{\partial J}{\partial x}.
\]
For linear state equations, the adjoint system is therefore simply the
transpose of the forward system.

\paragraph{Derivative with respect to the parameters}
For each component of the residual,
\[
F_i(x,p) = \sum\nolimits_{k=1}^{n} A_{ik}(p)x_k + b_i(p),
\]
and therefore
\[
\frac{\partial F_i}{\partial p_j} = \sum\nolimits_{k=1}^{n} \frac{\partial A_{ik}}{\partial p_j}x_k + \frac{\partial b_i}{\partial p_j}.
\]
In vector form,
\[
\frac{\partial F}{\partial p_j} = \frac{\partial A}{\partial p_j}x + \frac{\partial b}{\partial p_j}, \qquad \in\mathbb R^{n\times 1}.
\]
The full parameter Jacobian is assembled column-by-column as
\[
\frac{\partial F}{\partial p} = \begin{bmatrix} \frac{\partial A}{\partial p_1}x+\frac{\partial b}{\partial p_1} & \cdots & \frac{\partial A}{\partial p_m}x+\frac{\partial b}{\partial p_m} \end{bmatrix}.
\]
Substituting these expressions into
\[
\nabla_p J = \frac{\partial J}{\partial p} + \left(\frac{\partial F}{\partial p}\right)^T \lambda
\]
produces the reduced gradient required by the optimization algorithm.
The linear case illustrates an important point: the adjoint problem is
not an additional physical model. It is obtained directly from the
transpose of the discretized forward operator. As a result, every
choice made in constructing the forward discretization determines the
structure of the adjoint system.
%%%%%%%%%%%%%%%%%%%%%%%%%%%
%%%%%%%%%%%%%%%%%%%%%%%%%%%
\section{Gradient-Based Optimization with Adjoint Gradients}
\label{sec:solver}

The backward recursion of Section~\ref{sec:backward} computes the
reduced gradient $\nabla_p J(p_\ell)$ at a cost of one forward solve and one
adjoint solve, independent of the number of parameters. With this gradient
in hand, the PDE-constrained problem~\eqref{eq:opt} is reduced to an
unconstrained problem in $p$ alone,
\(
\min_{p} J(p),
\)
in which the governing equations enter only through the evaluation of
$J(p_\ell)$ and $\nabla_p J(p_\ell)$ at each iterate. The role of the
adjoint method is confined to this evaluation; it does not select the
optimization algorithm, and the remainder of this section treats that
algorithm as a separate, standard choice built on top of it.

\paragraph{Step length as an independent requirement} A gradient-based
method generates a sequence of iterates
\[
p_{\ell+1} = p_\ell + \alpha_\ell d_\ell,
\]
where $d_\ell$ is a descent direction satisfying $g_\ell^\top d_\ell < 0$,
with $g_\ell = \nabla_p J(p_\ell)$, and $\alpha_\ell > 0$ is a step length.
The descent condition established in Section~\ref{sec:adjoint} guarantees
only that some interval $(0,\bar\alpha_\ell]$ exists on which
$J(p_\ell+\alpha_\ell d_\ell) \le J(p_\ell)$; it is silent on the value of
$\bar\alpha_\ell$, which depends on the curvature of $J$ along $d_\ell$ and
is fixed by the line search, not by the adjoint computation that produced
$d_\ell$. In practice $\bar\alpha_\ell$ is never computed exactly; it is
instead approximated adaptively, by choosing $\alpha_\ell$ from the local
behavior of $J$ and $\nabla_p J$ at the current iterate rather than fixing
it in advance.\footnote{A constant step length is occasionally used for its simplicity,
but it requires an a priori bound on the curvature of $J$ to guarantee
descent at every iterate, a bound rarely available for a PDE-constrained
objective; chosen conservatively enough to be safe, it is typically far
smaller than an adaptive method would select, and chosen otherwise it
offers no guarantee at all. For an introduction to the standard
mechanisms for adaptive step-length selection---backtracking and
Wolfe-condition line search, trust-region methods, and their convergence
guarantees---see~\cite[Chs.~3--4]{nocedal} or~\cite[Ch.~1--2]{bertsekas2016}.}
For a PDE-constrained objective this adaptivity is more than a matter of
convergence rate: $J(p)$ is defined through the solution of $F(x,p)=0$, so
a step outside $(0,\bar\alpha_\ell]$ can move the trial parameter into a
region where the forward problem itself is unstable or fails to converge,
and the objective value the line search requires is then unavailable
rather than merely large. Section~\ref{sec:discrete-wave} exhibits this
failure mode directly, together with the two safeguards it necessitates in
practice: a feasibility check in the forward solver, and a
sufficient-decrease test applied to the realized, post-projection step
rather than to the nominal one. Correctness of $d_\ell$ and correctness of
the line search are independent requirements, and the adjoint method
establishes only the former.

\paragraph{The search direction as a separate choice} The simplest choice of $d_\ell$ is $d_\ell = -g_\ell$, giving gradient descent,
\[
p_{\ell+1} = p_\ell - \alpha_\ell \nabla_p J(p_\ell).
\]
Because it uses no curvature information, its convergence rate is linear
at best and is governed by the conditioning of $J$ near $p_\ell$, which
makes it a poor choice once $p$ has more than a few hundred components.
Newton's method recovers superlinear convergence by solving
$\nabla^2_p J(p_\ell)\, d_\ell = -g_\ell$, but forming or factoring the
Hessian $\nabla^2_p J(p_\ell) \in \mathbb{R}^{m\times m}$ is impractical
once $m$ reaches the thousands to millions typical of a discretized
parameter field. The limited-memory BFGS (L-BFGS) method occupies the
space between these two: it builds an implicit low-rank approximation
$B_\ell$ of the inverse Hessian from the $M$ most recent pairs
$(s_i, y_i) = (p_{i+1}-p_i,\, g_{i+1}-g_i)$, at a storage and computational
cost of $\mathcal{O}(Mm)$ rather than $\mathcal{O}(m^2)$, and it is the
method used in the demonstration of Section~\ref{sec:discrete-wave}.

At every iteration, L-BFGS interacts with the governing equations solely
through the gradient evaluation of Section~\ref{sec:backward};
everything else proceeds exactly as in an unconstrained optimization
problem with no reference to $F$. The complete procedure is summarized
below.

\begin{tcolorbox}[colback=gray!5,colframe=UCIblue,
  title=L-BFGS Optimization with Adjoint Gradients]
\begin{enumerate}[leftmargin=2em]
  \item \textbf{Initialization.} Choose an initial parameter estimate
  $p_0$, a memory size $M$, and a termination tolerance $\varepsilon$.
  Set $\ell=0$.
  \item \textbf{Gradient Evaluation.} Given $p_\ell$,
  \begin{enumerate}
    \item solve the forward problem $F(x,p_\ell)=0$ to obtain the state
    $x_\ell$;
    \item solve the adjoint equation
    $\displaystyle \left(\frac{\partial F}{\partial x}\right)^{\!\top}\!
    \lambda_\ell = -\frac{\partial J}{\partial x}$;
    \item compute the reduced gradient
    $\displaystyle g_\ell = \nabla_p J(p_\ell) = \frac{\partial J}{\partial p}
    + \left(\frac{\partial F}{\partial p}\right)^{\!\top}\! \lambda_\ell$.
  \end{enumerate}
  \item \textbf{Search Direction.} Compute $d_\ell = -B_\ell g_\ell$,
  where $B_\ell$ is the L-BFGS inverse-Hessian approximation obtained
  from the two-loop recursion below.
  \item \textbf{Line Search.} Determine a step length $\alpha_\ell$
  satisfying the Wolfe conditions
  \[
  J(p_\ell+\alpha_\ell d_\ell) \;\le\; J(p_\ell) + c_1 \alpha_\ell
  g_\ell^\top d_\ell,
  \qquad
  \nabla_p J(p_\ell+\alpha_\ell d_\ell)^\top d_\ell \;\ge\; c_2\, g_\ell^\top
  d_\ell.
  \]
  \item \textbf{Parameter Update.} Set
  $p_{\ell+1} = p_\ell + \alpha_\ell d_\ell$.
  \item \textbf{Memory Update.} Form the correction pair
  $s_\ell = p_{\ell+1}-p_\ell$, $y_\ell = g_{\ell+1}-g_\ell$; store
  $(s_\ell,y_\ell)$ and discard the oldest pair once more than $M$ are
  retained.
  \item \textbf{Termination.} If $\|\nabla_p J(p_\ell)\| < \varepsilon$,
  stop; otherwise set $\ell \leftarrow \ell+1$ and return to Step 2.
\end{enumerate}
\end{tcolorbox}

The only step at which the governing equations are consulted is Step 2;
Steps 3--6 depend on $F$ only through the values $g_\ell$ it has already
produced. The search direction of Step 3 is obtained from the two-loop
recursion,

\begin{tcolorbox}[colback=gray!5,colframe=UCIblue,
  title=L-BFGS Two-Loop Recursion]
Given $g_\ell$ and the stored pairs $\{(s_i,y_i)\}_{i=\ell-M}^{\ell-1}$,
set $q \leftarrow g_\ell$.

\smallskip
\noindent\textit{Backward pass}, for $i = \ell-1, \ell-2, \dots, \ell-M$:
\[
\rho_i \leftarrow (y_i^\top s_i)^{-1},
\qquad
\alpha_i \leftarrow \rho_i\, s_i^\top q,
\qquad
q \leftarrow q - \alpha_i\, y_i.
\]

\smallskip
\noindent Set
\[
\gamma \leftarrow \frac{s_{\ell-1}^\top y_{\ell-1}}{y_{\ell-1}^\top y_{\ell-1}},
\qquad
r \leftarrow \gamma\, q.
\]

\smallskip
\noindent\textit{Forward pass}, for $i = \ell-M, \ell-M+1, \dots, \ell-1$:
\[
\beta \leftarrow \rho_i\, y_i^\top r,
\qquad
r \leftarrow r + s_i(\alpha_i - \beta).
\]

\smallskip
\noindent Return $d_\ell = -r$.
\end{tcolorbox}

Retaining only the $M$ most recent pairs bounds the memory and per-iteration
cost at $\mathcal{O}(Mm)$, in contrast to the $\mathcal{O}(m^2)$ storage a
dense quasi-Newton update would require; this is what makes L-BFGS
suitable for a parameter dimension $m$ ranging from thousands to millions,
a regime in which forming $\nabla_p^2 J$ even once is already out of reach.
% ========================================================
% ========================================================
% ========================================================
\section{Discrete Wave Equation and Its Adjoint: A Step-by-Step Construction}\label{sec:discrete-wave}

In this section, we construct the adjoint method explicitly for the vibrating
string introduced in Section~\ref{sec:inverse}, working entirely at the
discrete level. This small system allows us to see exactly how the forward
model evolves, how discretization choices change the algebraic structure, and
how the adjoint propagates sensitivities backward through time.

Recall from Section~\ref{sec:inverse} that the governing equation is the
one-dimensional wave equation in divergence form, \eqref{eq:pde}, and that the
parameter $p$ consists of the discrete values of $c(z)^2$.

To express the wave equation as a system suitable for time stepping, we first
rewrite it as a first-order system in time:
\[
\dot{u} = v, \qquad
\dot{v} = \mathbf{L}(c^2)\,u,
\]
where $\mathbf{L}(c^2)$ denotes the spatial operator
\[
\mathbf{L}(c^2)u \;=\; \frac{\partial}{\partial z}\left(c^2 \frac{\partial u}{\partial z}\right),
\]
still acting on the continuous function $u(z,t)$; conceptually, it represents
how spatial variations in the displacement $u$ generate acceleration through
the wave speed $c(z)$. Spatial discretization, carried out next, replaces $u$
by its values on a grid and $\mathbf{L}(c^2)$ by a matrix that approximates
this operator, with the approximation becoming exact as the grid is refined.

\paragraph{Initial Conditions}
The spatial boundary conditions are the Dirichlet conditions already
prescribed in Section~\ref{sec:inverse}: the string is fixed at both ends,
$u(0,t) = u(L,t) = 0$. To carry the construction through concretely, we also fix a specific
initial state: the string is plucked to a height $h$ at an interior
point $z=L_1$, with $0 < L_1 < L$, and released from rest, so that
\begin{equation}\label{eq:initial_disp}
u(z,0) =
\begin{cases}
\dfrac{h}{L_1}\,z, & 0 \le z \le L_1, \\[4pt]
\dfrac{h}{L-L_1}\,(L-z), & L_1 < z \le L,
\end{cases}
\qquad v(z,0) = 0.
\end{equation}
As noted in Section~\ref{sec:inverse}, this initial state is given data, not
an optimization variable: it is the starting point of the forward simulation
and influences every later state, but it has no associated adjoint variable
in the backward computation. Its discrete counterpart, the initial state
vector $x^0$, is given once the spatial grid is introduced below.

%=========================================================
\subsection{Collocated Discretization (State Variables on the Same Grid)}\label{sec:collocated}
This subsection carries out the first of the two discretizations mentioned above: all state variables---displacement and
velocity alike---are placed on the same spatial and temporal grid, and the forward model, the state equation
$F(x,p)=0$, and the resulting adjoint recursion are constructed explicitly
for this choice. Rather than starting from a large abstract system, we
build the model step by step for a small example, with the grid
deliberately small---and, to keep spatial and temporal indices visually
distinct throughout, of a different size in each direction---so that every
vector and matrix appearing below can be written out in full:
\[
n_z = 7 \quad \text{spatial nodes}, \qquad n_t = 4 \quad \text{time levels}.
\]
Here $n_z$ counts the spatial nodes $z_0,\dots,z_6$ inclusive of the two
boundary points, and $n_t$ counts the time levels $t^0,\dots,t^3$ inclusive
of the prescribed initial level; the domain therefore has five interior
nodes $z_1,\dots,z_5$, spaced $\Delta z = L/(n_z-1)$ apart, and three
computed time steps, reaching a terminal level $N=3$. The boundary values
inherited from Section~\ref{sec:inverse} read, at every time level,
\[
u_0^k = 0, \qquad u_6^k = 0, \qquad
v_0^k = 0, \qquad v_6^k = 0, \qquad k=1,2,3,
\]
the velocity conditions following from the displacement conditions: since
$u(0,t)=u(L,t)=0$ for all $t$, the boundary points do not move at all, so
their velocity vanishes identically as well. The interior values at the
initial level, $u_i^0$, are specified by the plucked profile
in~\eqref{eq:initial_disp}, and the string is released from rest, so
\[
u_i^0 = u(z_i,0), \qquad v_i^0 = 0, \qquad i=0,\dots,6.
\]

At each time level the state collects the interior displacements and
velocities,
\begin{equation}
x^k =
\begin{bmatrix}
u_1^k & u_2^k & u_3^k & u_4^k & u_5^k &
v_1^k & v_2^k & v_3^k & v_4^k & v_5^k
\end{bmatrix}^\top
\in \mathbb{R}^{10}, \qquad k=1,2,3.
\end{equation}
The full space-time trajectory is thus the collection $x^1,x^2,x^3$,
together with the given initial state
\[
x^0=
\begin{bmatrix}
u_1^0 & u_2^0 & u_3^0 & u_4^0 & u_5^0 &
0 & 0 & 0 & 0 & 0
\end{bmatrix}^\top.
\]

\paragraph{The discrete spatial operator}
Approximating the continuous operator $\mathbf{L}(c^2)$ of
Section~\ref{sec:discrete-wave} by centered differences gives, at any node $i$,
\begin{equation}\label{eq:L-stencil}
(\mathbf{L}(c^2)u)_i \;\approx\;
\frac{1}{\Delta z^2}
\left[
\left(\frac{c^2_{i+1}+c^2_i}{2}\right)(u_{i+1}-u_i)
-
\left(\frac{c^2_i+c^2_{i-1}}{2}\right)(u_i-u_{i-1})
\right],
\end{equation}
which involves the two spatial neighbors $u_{i-1}$ and $u_{i+1}$ together
with $u_i$ itself. For $i=2,3,4$ both neighbors are interior unknowns, and
\eqref{eq:L-stencil} contributes only to the $5\times5$ interior-to-interior
block of the operator. At $i=1$ and $i=5$, however, one of the two
neighbors is a boundary node: writing~\eqref{eq:L-stencil} out for $i=1$,
\begin{equation}
(\mathbf{L}(c^2)u)_1 =
\frac{1}{2\Delta z^2}\Big[-(c_0^2+2c_1^2+c_2^2)\,u_1 + (c_1^2+c_2^2)\,u_2\Big]
\;+\; \frac{c_0^2+c_1^2}{2\Delta z^2}\,u_0,
\end{equation}
and symmetrically for $i=5$ with $u_6$ in place of $u_0$. The last term in
each case is precisely the coupling to the boundary that a node-by-node
reading of~\eqref{eq:L-stencil} must not lose. Because
Section~\ref{sec:inverse} prescribes $u_0^k=u_6^k=0$ for every $k$, these
two terms vanish identically---not because they are dropped, but because
the boundary data driving them is zero---and what remains is the
$5\times5$ tridiagonal matrix acting on the interior vector $u^k$ alone,
\begin{equation}\label{eq:L-collocated}
\mathbf{L}(c^2) = \frac{1}{2\Delta z^2}
\begin{bmatrix}
-(c^2_0+2c^2_1+c^2_2) & c^2_1+c^2_2 & 0 & 0 & 0 \\
c^2_1+c^2_2 & \!\!-(c^2_1+2c^2_2+c^2_3) & c^2_2+c^2_3 & 0 & 0 \\
0 & c^2_2+c^2_3 & \!\!-(c^2_2+2c^2_3+c^2_4) & c^2_3+c^2_4 & 0 \\
0 & 0 & c^2_3+c^2_4 & \!\!-(c^2_3+2c^2_4+c^2_5) & c^2_4+c^2_5 \\
0 & 0 & 0 & c^2_4+c^2_5 & \!\!-(c^2_4+2c^2_5+c^2_6)
\end{bmatrix}. 
\end{equation}
The boundary values $c_0^2,c_6^2$ still appear, in the diagonal entries for
$i=1$ and $i=5$: they are wave-speed data, needed to average the material
property across the boundary interval, and are retained even though the
boundary \emph{displacements} $u_0,u_6$ that would otherwise multiply them
are zero. Since the inverse problem of Section~\ref{sec:inverse} treats
only the interior profile as unknown, $c_0^2$ and $c_6^2$ are prescribed
data rather than components of the parameter vector
$p=[c_1^2,\dots,c_5^2]^\top\in\mathbb{R}^5$, which appears through the
interior tridiagonal coefficients of $\mathbf{L}(c^2)$ alone. Note also
that $\mathbf{L}(c^2)$, as constructed, is symmetric; this will simplify
the adjoint recursion below. Finally, neither update below ever
references a neighboring \emph{velocity}: $v_0^k$ and $v_6^k$ do not enter
$\mathbf{L}(c^2)$ or the state equation at all, prescribed or not.

\paragraph{Time stepping and the state equation}
Advancing $u$ and $v$ by forward Euler for the collocated discretization
gives, componentwise, for $i=1,\dots,5$,
\begin{subequations}
\begin{align}
u_i^{k+1} &= u_i^k + \Delta t\, v_i^k, \\
v_i^{k+1} &= v_i^k + \Delta t\, \big(\mathbf{L}(c^2)u^k\big)_i.
\end{align}
\end{subequations}
The first equation is purely local: $u_i^{k+1}$ depends only on $u_i^k$ and
$v_i^k$, at the same node, with no spatial coupling at all. The second
equation is where the spatial coupling derived above enters: $v_1^{k+1}$
depends on $u_0^k,u_1^k,u_2^k$, and $v_5^{k+1}$ depends on
$u_4^k,u_5^k,u_6^k$, with the $u_0^k$ and $u_6^k$ contributions equal to
zero for the reason given above. Consequently, with $\mathbf{L}(c^2)$ as
in~\eqref{eq:L-collocated}, the two equations stack \emph{exactly}, with no
omitted term, into
\begin{equation}\label{eq:A-collocated}
x^{k+1} = A(p)\,x^k, \qquad
A(p) =
\begin{bmatrix}
\mathbf{I} & \Delta t\,\mathbf{I} \\
\Delta t\,\mathbf{L}(c^2) & \mathbf{I}
\end{bmatrix}.
\end{equation}
For the three time steps of this example, the forward evolution reads
$x^1=A(p)x^0$, $x^2=A(p)x^1$, $x^3=A(p)x^2$.

\paragraph{The full space-time system}
Collecting all three time levels into a single system expresses the
constraint $F(x,p)=0$ explicitly:
\begin{equation}\label{eq:F-collocated}
\begin{bmatrix}
I & 0 & 0 \\
-A(p) & I & 0 \\
0 & -A(p) & I
\end{bmatrix}
\begin{bmatrix} x^1 \\ x^2 \\ x^3 \end{bmatrix}
=
\begin{bmatrix} A(p)x^0 \\ 0 \\ 0 \end{bmatrix}.
\end{equation}
The first block row enforces the initial step $x^1=A(p)x^0$; because $x^0$
is known, it enters as data on the right-hand side rather than as an
unknown, consistent with Section~\ref{sec:discrete-wave}. The remaining rows
enforce the recursive time-marching structure. The matrix on the left is
block lower triangular, so its dependence on $x^0$ propagates forward
through the entire trajectory, exactly the causal structure noted in the
introduction.

\paragraph{The adjoint equation}
Specializing the general adjoint equation of Section~\ref{sec:backward},
$(\partial F/\partial x)^\top\lambda = -\partial J/\partial x$, to the
block lower-triangular constraint~\eqref{eq:F-collocated} gives a block
upper-triangular system,
\begin{equation}\label{eq:adjoint-collocated}
\begin{bmatrix}
I & -A(p)^\top & 0 \\
0 & I & -A(p)^\top \\
0 & 0 & I
\end{bmatrix}
\begin{bmatrix} \lambda^1 \\ \lambda^2 \\ \lambda^3 \end{bmatrix}
=
-\begin{bmatrix}
\nabla_{x^1}J \\ \nabla_{x^2}J \\ \nabla_{x^3}J
\end{bmatrix},
\end{equation}
which, as the introduction anticipated, is solved by back substitution:
from the terminal condition
\begin{equation}
\lambda^3 = -\nabla_{x^3}J,
\end{equation}
one proceeds by the backward recursion
\begin{align}
\lambda^2 &= A(p)^\top\lambda^3 - \nabla_{x^2}J, \\
\lambda^1 &= A(p)^\top\lambda^2 - \nabla_{x^1}J.
\end{align}

\paragraph{Backward recursion in components}
The recursion is initialized at the terminal time level $k=N$ (here $N=3$)
with
\[
\lambda^N = -\nabla_{x^N}J\Big|_{u^N,v^N}\,\,,
\]
the last block row of~\eqref{eq:adjoint-collocated}, since there is no
level $N+1$ to supply an incoming term $A(p)^\top\lambda^{N+1}$: the row
reduces to $I\lambda^N=-\nabla_{x^N}J$ with nothing else to solve for.
Written this way, the terminal condition holds for any differentiable
objective $J$, evaluated at the state $x^N=(u^N,v^N)$ that the forward
simulation has already produced; unlike every earlier $\lambda^k$, no
recursion is needed to obtain it.

For the quadratic misfit of Section~\ref{sec:inverse},
$J(x,p)=\tfrac12\|Hx-d\|^2+\tfrac{\beta}{2}\mathcal{R}(p)$, this general
expression can be opened up further. Writing $H_N$ for the rows of $H$
that select the components of $x^N$ observed at the final recording
instant and $d_N$ for the corresponding measured values,
\[
\lambda^N = -\nabla_{x^N}J\Big|_{u^N,v^N} = -H_N^\top\big(H_N x^N - d_N\big),
\]
a vector that is zero in every component of $x^N$ not observed by a
sensor. This is the same mechanism that fixed the boundary values
$u_0^k=u_6^k=v_0^k=v_6^k=0$ in the forward problem: there, a component of
the state was pinned to zero because it was excluded from the interior
unknowns by the Dirichlet condition; here, a component of the terminal
adjoint is pinned to zero because it is excluded from the data term by
$H$. Once the forward simulation has produced $x^N$, computing $\lambda^N$
requires only this evaluation, quadratic cost or not.

Writing each $\lambda^k$ in its displacement and velocity blocks,
$\lambda^k = (\lambda_u^k,\lambda_v^k)$, and using
$A(p)^\top = \begin{bmatrix} I & \Delta t\,\mathbf{L}(c^2)^\top \\ \Delta t\,I & I\end{bmatrix}$,
the recursion for $k=N-1,\dots,1$ becomes
\begin{subequations}
\begin{align}
\lambda_u^k &= \lambda_u^{k+1} + \Delta t\,\mathbf{L}(c^2)^\top\lambda_v^{k+1}
- \frac{\partial J}{\partial u^k}\Big|_{u^k,v^k}, \\
\lambda_v^k &= \lambda_v^{k+1} + \Delta t\,\lambda_u^{k+1}
- \frac{\partial J}{\partial v^k}\Big|_{u^k,v^k},
\end{align}
\end{subequations}
and, since $\mathbf{L}(c^2)$ is symmetric, $\mathbf{L}(c^2)^\top=\mathbf{L}(c^2)$
in both expressions. These two equations make explicit how sensitivity
information flows backward: the velocity adjoint at step $k$ is driven by
the displacement adjoint one step ahead, and vice versa, mirroring---in
reverse---the coupling between $u$ and $v$ in the forward recursion.

\paragraph{The role of the initial and boundary conditions}
As in Section~\ref{sec:discrete-wave}, there is no adjoint variable associated
with $x^0$, since $x^0$ is prescribed data rather than an unknown. Its
influence on the objective is nonetheless present: it enters the forward
trajectory through $x^1=A(p)x^0$, and the adjoint variable $\lambda^1$
carries the corresponding sensitivity backward. The initial condition thus
affects the computed gradient only indirectly, through the state it
generates and the adjoint that propagates from it.

The boundary values $u_0^k,v_0^k,u_6^k,v_6^k$ play an analogous role.
They too are excluded from $x^k$ and never acquire adjoint variables of
their own, for the same reason $x^0$ does not: they are prescribed data,
not unknowns solved for by $F(x,p)=0$. Their influence on the computed
gradient is threaded through $\mathbf{L}(c^2)$ instead of through an
explicit adjoint component: the boundary displacements $u_0^k,u_6^k$ were
shown in the derivation of $\mathbf{L}(c^2)$ to enter the interior update
only through the terms that vanish because $u_0^k=u_6^k=0$, while the
boundary wave speeds $c_0^2,c_6^2$ remain, as data, in the diagonal
entries of $\mathbf{L}(c^2)$ itself. Consequently the boundary conditions
affect $\lambda^k$ only through the matrix $\mathbf{L}(c^2)^\top$ that
appears in every step of the backward recursion, never through a
boundary adjoint term that would need to be added or initialized
separately---exactly as no such term needed to be added when the forward
state equation was assembled.

\paragraph{The parameter gradient}
Only the $(2,1)$ block of $A(p)$ in~\eqref{eq:A-collocated} depends on
$p$, so for $p_j=c_j^2$,
\[
\frac{\partial A}{\partial p_j} =
\begin{bmatrix}
0 & 0 \\[2pt]
\Delta t\,\dfrac{\partial \mathbf{L}(c^2)}{\partial p_j} & 0
\end{bmatrix},
\]
with $\partial \mathbf{L}(c^2)/\partial p_j$ obtained by
differentiating~\eqref{eq:L-collocated} entry by entry: each interior
parameter $c_j^2$ enters through the diagonal of row $j$ and, through
the face-averaging in~\eqref{eq:L-stencil}, the diagonals of its two
immediate neighbors and the off-diagonal entries connecting them---three
diagonal entries and two symmetric off-diagonal pairs for $j=2,3,4$,
fewer at the domain edges $j=1,5$, where one neighbor is the
prescribed, parameter-free boundary node. Writing $\lambda^{k+1}$ in its
displacement and velocity blocks as before,
\[
\big(\lambda^{k+1}\big)^\top\frac{\partial A}{\partial p_j}x^k
= \Delta t\,\big(\lambda_v^{k+1}\big)^\top
\frac{\partial \mathbf{L}(c^2)}{\partial p_j}\,u^k,
\]
and the reduced gradient is assembled, one time step at a time, as
\begin{equation}\label{eq:grad-collocated}
(\nabla_p J)_j = \frac{\partial J}{\partial p_j}
- \Delta t\sum\nolimits_{k=0}^{N-1}\big(\lambda_v^{k+1}\big)^\top
\frac{\partial \mathbf{L}(c^2)}{\partial p_j}\,u^k.
\end{equation}
This is a zero-lag cross-correlation between the forward displacement
field $u^k$ and the adjoint velocity field $\lambda_v^{k+1}$,
accumulated over the whole record---the discrete \emph{imaging
condition} that full-waveform-inversion codes assemble to update a
velocity model from a forward and an adjoint wavefield.
%%%%%%%%%%%%%%%%%%%%%%%%%%
%%%%%%%%%%%%%%%%%%%%%%%%%%%
% ========================================================

%%%%%%%%%%%%%%%%%
%%%%%%%%%%%%%%%%%
\subsection{Staggered Discretization (State Variables on Offset Grids)}\label{sec:staggered}
We now carry out the second of the two discretizations methods we mentioned earlier. Rather than
placing displacement and velocity on the same grid, as in
Section~\ref{sec:collocated}, the staggered scheme offsets velocity by half
a time step and the wave-speed coefficient by half a spatial step, so that
each derivative is approximated at the point where it is naturally defined.
This is a common device for reducing phase error and spurious numerical
dissipation in wave problems, and it is adopted here not for that benefit
but because it changes the bookkeeping of the adjoint recursion in a way
the collocated scheme does not exhibit---the physics and the boundary data
are exactly as in Section~\ref{sec:inverse}; only the discrete algorithm
changes.

We use the same grid as in Section~\ref{sec:collocated}: $n_z=7$ spatial
nodes $z_0,\dots,z_6$ inclusive of the two boundary points, spaced
$\Delta z=L/(n_z-1)$ apart, and $n_t=4$ time levels $t^0,\dots,t^3$
inclusive of the prescribed initial level, giving three computed time
steps and a terminal level $N=3$. Displacement $u$ is defined at the
integer nodes $z_1,\dots,z_5$ and integer levels $t^1,t^2,t^3$, exactly as
in the collocated case; velocity $v$ is instead defined at the half time
levels $t^{1/2},t^{3/2},t^{5/2}$; and the coefficient $c^2$ is defined at
the spatial half-nodes $z_{1/2},z_{3/2},\dots,z_{11/2}$, one between each
pair of neighboring integer nodes. The boundary and initial data are
unchanged from Section~\ref{sec:collocated}:
\[
u_0^k = 0, \qquad u_6^k = 0, \qquad k=1,2,3,
\]
with the corresponding boundary velocities $v_0^{k+1/2}=v_6^{k+1/2}=0$ for
$k=0,1,2$, since the boundary points do not move at any time, integer or
half-integer. The interior displacements at the initial level are again
$u_i^0=u(z_i,0)$ from the plucked profile~\eqref{eq:initial_disp}, and
the string is released from rest; the leapfrog recursion below is
started with $v_i^{-1/2}=0$, the simplest choice consistent with
$v(z,0)=0$. This start is only first-order accurate at the first
step---the second-order-consistent choice is the half-step Taylor
correction $v_i^{-1/2}=-\tfrac{\Delta t}{2}\big(\mathbf
L_{\mathrm{stag}}(c^2)u^0\big)_i$, and production leapfrog codes use
it---but the choice of start-up has no bearing on this tutorial's
argument: whichever forward recursion is coded, and however accurately
it approximates the continuous wave equation, the adjoint is the exact
transpose of \emph{that recursion}, not of the equation it
approximates. The simpler start is used here only to keep the worked
example free of an extra half-step formula; because it perturbs a
single step by $O(\Delta t^2)$, no larger than the per-step local error
of a second-order scheme already carries, it does not change the order
of accuracy of the simulation used in
Section~\ref{sec:numerical-demonstration}.

With this structure, the forward update is
\begin{subequations}
\begin{align}
v^{k+1/2} &= v^{k-1/2} + \Delta t\,\mathbf{L}_{\mathrm{stag}}(c^2)\,u^k, \label{eq:stag-v}\\
u^{k+1} &= u^k + \Delta t\, v^{k+1/2}. \label{eq:stag-u}
\end{align}
\end{subequations}
Unlike the collocated scheme, advancing one time level now requires two
sequential sub-steps: velocity is updated first from $u^k$, and displacement
is then updated from the newly computed $v^{k+1/2}$,
\[
u^k \;\longrightarrow\; v^{k+1/2} \;\longrightarrow\; u^{k+1}.
\]
At each time level the state collects the displacement and the velocity
one half-step behind it,
\begin{equation}
x^k =
\begin{bmatrix}
u_1^k & u_2^k & u_3^k & u_4^k & u_5^k &
v_1^{k-1/2} & v_2^{k-1/2} & v_3^{k-1/2} & v_4^{k-1/2} & v_5^{k-1/2}
\end{bmatrix}^\top
\in \mathbb{R}^{10}, \qquad k=1,2,3,
\end{equation}
and the given initial state is
\[
x^0=
\begin{bmatrix}
u_1^0 & u_2^0 & u_3^0 & u_4^0 & u_5^0 &
0 & 0 & 0 & 0 & 0
\end{bmatrix}^\top.
\]

\paragraph{The discrete spatial operator}
Component-wise,
\begin{equation}\label{eq:L-stag-stencil}
(\mathbf{L}_{\mathrm{stag}}(c^2)u)_i =
\frac{c_{i+1/2}^2(u_{i+1}-u_i) - c_{i-1/2}^2(u_i-u_{i-1})}{\Delta z^2},
\end{equation}
again a discrete representation of $\partial_z(c^2\partial_z u)$, this time
built directly from wave speeds sampled at the half-nodes, with no
averaging needed since $c^2$ already lives exactly where the flux is
evaluated. As in Section~\ref{sec:collocated}, this involves the two
spatial neighbors $u_{i-1},u_{i+1}$; for $i=1$ one of them is the boundary
value $u_0$. Writing~\eqref{eq:L-stag-stencil} out for $i=1$,
\begin{equation}
(\mathbf{L}_{\mathrm{stag}}(c^2)u)_1 =
\frac{1}{\Delta z^2}\Big[-(c_{1/2}^2+c_{3/2}^2)\,u_1 + c_{3/2}^2\,u_2\Big]
\;+\; \frac{c_{1/2}^2}{\Delta z^2}\,u_0,
\end{equation}
and symmetrically for $i=5$ with $u_6$ and $c_{11/2}^2$ in place of $u_0$
and $c_{1/2}^2$. As in the collocated case, these last terms vanish
identically because $u_0^k=u_6^k=0$, and what remains is the $5\times5$
tridiagonal matrix acting on the interior vector $u^k$ alone,
\begin{equation}\label{eq:L-stag}
\mathbf{L}_{\mathrm{stag}}(c^2) = \frac{1}{\Delta z^2}
\begin{bmatrix}
-(c^2_{1/2}+c^2_{3/2}) & c^2_{3/2} & 0 & 0 & 0 \\
c^2_{3/2} & \!\!-(c^2_{3/2}+c^2_{5/2}) & c^2_{5/2} & 0 & 0 \\
0 & c^2_{5/2} & \!\!-(c^2_{5/2}+c^2_{7/2}) & c^2_{7/2} & 0 \\
0 & 0 & c^2_{7/2} & \!\!-(c^2_{7/2}+c^2_{9/2}) & c^2_{9/2} \\
0 & 0 & 0 & c^2_{9/2} & \!\!-(c^2_{9/2}+c^2_{11/2})
\end{bmatrix},
\end{equation}
symmetric, exactly as $\mathbf{L}(c^2)$ was in Section~\ref{sec:collocated}.
Unlike that case, no boundary wave speed is needed here at all: every
coefficient in~\eqref{eq:L-stag} sits strictly inside the domain, at a
half-node $z_{1/2},\dots,z_{11/2}$, so all six values
\begin{equation}
p = \big[\,c^2_{1/2},\,c^2_{3/2},\,c^2_{5/2},\,c^2_{7/2},\,c^2_{9/2},\,c^2_{11/2}\,\big]^\top
\end{equation}
are unknowns---one more component than the five interior values
$p=[c_1^2,\dots,c_5^2]^\top$ of the collocated parameterization, and none
of them prescribed boundary data, precisely the second choice anticipated
in Section~\ref{sec:inverse}. As before, neither update below ever
references a neighboring \emph{velocity}: $v_0^{k+1/2}$ and $v_6^{k+1/2}$
do not enter $\mathbf{L}_{\mathrm{stag}}(c^2)$ or the state equation at
all.

\paragraph{Time stepping and the state equation}
Componentwise, for $i=1,\dots,5$,
\begin{subequations}
\begin{align}
v_i^{k+1/2} &= v_i^{k-1/2} + \Delta t\,\big(\mathbf{L}_{\mathrm{stag}}(c^2)u^k\big)_i, \\
u_i^{k+1} &= u_i^k + \Delta t\, v_i^{k+1/2}.
\end{align}
\end{subequations}
The velocity update is where the spatial coupling derived above enters:
$v_1^{k+1/2}$ depends on $u_0^k,u_1^k,u_2^k$, and $v_5^{k+1/2}$ depends on
$u_4^k,u_5^k,u_6^k$, with the $u_0^k$ and $u_6^k$ contributions equal to
zero for the reason given above. The displacement update, by contrast, is
purely local: $u_i^{k+1}$ depends only on $u_i^k$ and the newly computed
$v_i^{k+1/2}$, at the same node. Eliminating the half-step velocity by
substituting the first equation into the second gives the single linear
map
\begin{equation}\label{eq:A-stag}
x^{k+1} = A(p)\,x^k, \qquad
A(p) = \begin{bmatrix}
I+\Delta t^2\,\mathbf{L}_{\mathrm{stag}}(c^2) & \Delta t\,I \\
\Delta t\,\mathbf{L}_{\mathrm{stag}}(c^2) & I
\end{bmatrix}.
\end{equation}
For the three time steps of this example, the forward evolution reads
$x^1=A(p)x^0$, $x^2=A(p)x^1$, $x^3=A(p)x^2$.

\paragraph{The full space-time system}
Collecting all three time levels into a single system gives exactly the
block lower-triangular constraint of Section~\ref{sec:collocated},
\begin{equation}\label{eq:F-stag}
\begin{bmatrix}
I & 0 & 0 \\
-A(p) & I & 0 \\
0 & -A(p) & I
\end{bmatrix}
\begin{bmatrix} x^1 \\ x^2 \\ x^3 \end{bmatrix}
=
\begin{bmatrix} A(p)x^0 \\ 0 \\ 0 \end{bmatrix},
\end{equation}
with $A(p)$ replaced by its staggered counterpart~\eqref{eq:A-stag}: the
two discretizations differ in the entries of $A(p)$, not in the shape of
$F(x,p)=0$.

\paragraph{The adjoint equation}
Specializing the general adjoint equation of Section~\ref{sec:backward},
$(\partial F/\partial x)^\top\lambda=-\partial J/\partial x$, to the block
lower-triangular constraint~\eqref{eq:F-stag} gives the same block
upper-triangular system as in Section~\ref{sec:collocated}, now with
$A(p)^\top$ built from~\eqref{eq:A-stag}: from the terminal condition
\begin{equation}
\lambda^3 = -\nabla_{x^3}J,
\end{equation}
one proceeds by the backward recursion
\begin{align}
\lambda^2 &= A(p)^\top\lambda^3 - \nabla_{x^2}J, \\
\lambda^1 &= A(p)^\top\lambda^2 - \nabla_{x^1}J.
\end{align}

\paragraph{Backward recursion in components}
As in Section~\ref{sec:collocated}, the terminal condition
\[
\lambda^N = -\nabla_{x^N}J\Big|_{u^N,v^{N-1/2}}
\]
is known as soon as the forward simulation has produced $x^N=(u^N,v^{N-1/2})$,
with $N=3$: it is the last block row of the adjoint system, with no
incoming term to solve for, and it holds for any differentiable $J$. For
the quadratic misfit of Section~\ref{sec:inverse}, writing $H_N$ for the
rows of $H$ that select the observed components of $x^N$ and $d_N$ for the
corresponding data,
\[
\lambda^N = -H_N^\top\big(H_N x^N - d_N\big),
\]
zero in every unobserved component, exactly as in the collocated case.

Writing each $\lambda^k$ in its displacement and velocity blocks,
$\lambda^k=(\lambda_u^k,\lambda_v^{k-1/2})$, and using
$A(p)^\top = \begin{bmatrix} I+\Delta t^2\,\mathbf{L}_{\mathrm{stag}}(c^2)^\top & \Delta t\,\mathbf{L}_{\mathrm{stag}}(c^2)^\top \\ \Delta t\,I & I \end{bmatrix}$,
the recursion for $k=N-1,\dots,1$ becomes
\begin{align}
\lambda_u^k &= \big(I+\Delta t^2\,\mathbf{L}_{\mathrm{stag}}(c^2)\big)\lambda_u^{k+1}
+ \Delta t\,\mathbf{L}_{\mathrm{stag}}(c^2)\,\lambda_v^{k+1/2}
- \frac{\partial J}{\partial u^k}\Big|_{u^k,v^{k-1/2}}, \label{eq:lam-u-stag}\\
\lambda_v^{k-1/2} &= \Delta t\,\lambda_u^{k+1} + \lambda_v^{k+1/2}
- \frac{\partial J}{\partial v^{k-1/2}}\Big|_{u^k,v^{k-1/2}}, \label{eq:lam-v-stag}
\end{align}
using the symmetry of $\mathbf{L}_{\mathrm{stag}}(c^2)$ to replace its
transpose by itself, exactly as in Section~\ref{sec:collocated}.

Equations~\eqref{eq:lam-u-stag}--\eqref{eq:lam-v-stag} are written in
eliminated form, exactly parallel to the eliminated forward
map~\eqref{eq:A-stag}, so both appear to draw on $\lambda^{k+1}$
directly and independently of one another. They are not independent:
substituting~\eqref{eq:lam-v-stag} into
\[
\lambda_u^k = \lambda_u^{k+1}
+ \Delta t\,\mathbf{L}_{\mathrm{stag}}(c^2)\,\lambda_v^{k-1/2}
- \frac{\partial J}{\partial u^k}\Big|_{u^k,v^{k-1/2}}
\]
in place of $\lambda_v^{k-1/2}$ reproduces~\eqref{eq:lam-u-stag}
exactly. Equation~\eqref{eq:lam-v-stag} is therefore already the
adjoint of the displacement sub-step~\eqref{eq:stag-u}, and the pair
above is the adjoint of the velocity sub-step~\eqref{eq:stag-v}
followed by the adjoint of the displacement sub-step---the forward
sub-steps $u^k\to v^{k+1/2}\to u^{k+1}$ retraced by the backward sweep
as $\lambda^{k+1}\to\lambda_v^{k-1/2}\to\lambda_u^k$, the same two
half-steps in reverse order. Writing $A=A_2A_1$ for the sub-step
decomposition of~\eqref{eq:A-stag}, with $A_1$ the velocity half-step
and $A_2$ the displacement step, this is simply $A^\top =
A_1^\top A_2^\top$: the adjoint of a composition is the composition of
the adjoints, in reverse order. What~\eqref{eq:lam-u-stag}--\eqref{eq:lam-v-stag}
obscures is not the ordering itself but the fact that it has already
been eliminated, in parallel with the forward map---one more
manifestation of the tutorial's central point: the discretization, not
the physics, decides how the adjoint recursion is organized. It is
also, concretely, why full-waveform-inversion codes compute the
adjoint wavefield by rerunning the same leapfrog solver backward in
time with the data misfit as a source term~\cite{tarantola1984,virieux2009,fichtner2011}: the adjoint of
the staggered leapfrog is the staggered leapfrog, run in reverse.

\paragraph{The role of the initial and boundary conditions}
As in Section~\ref{sec:collocated}, there is no adjoint variable associated
with $x^0$; its influence on the objective is carried entirely through the
forward trajectory it generates and the backward recursion that returns
from it. The boundary values play the same role here as in the collocated
case, and by an even simpler mechanism: the boundary displacements
$u_0^k,u_6^k$ enter $\mathbf{L}_{\mathrm{stag}}(c^2)$ only through the terms
shown above to vanish because $u_0^k=u_6^k=0$, and---unlike the collocated
operator, whose diagonal retained the boundary wave speeds $c_0^2,c_6^2$ as
fixed data---no boundary material property enters
$\mathbf{L}_{\mathrm{stag}}(c^2)$ at all, since every coefficient it uses
sits strictly inside the domain. The boundary conditions therefore affect
$\lambda^k$ only through $\mathbf{L}_{\mathrm{stag}}(c^2)^\top$ in the
recursion above, never through a boundary adjoint term of their own,
exactly as in Section~\ref{sec:collocated}.

\paragraph{The parameter gradient}
Only the two blocks of $A(p)$ built from $\mathbf{L}_{\mathrm{stag}}(c^2)$
in~\eqref{eq:A-stag} depend on $p$, so for $p_j=c^2_{j+1/2}$,
\[
\frac{\partial A}{\partial p_j} =
\begin{bmatrix}
\Delta t^2\,\dfrac{\partial \mathbf{L}_{\mathrm{stag}}}{\partial p_j} & 0 \\[4pt]
\Delta t\,\dfrac{\partial \mathbf{L}_{\mathrm{stag}}}{\partial p_j} & 0
\end{bmatrix},
\]
with $\partial\mathbf{L}_{\mathrm{stag}}/\partial p_j$ the four-entry
matrix obtained by differentiating~\eqref{eq:L-stag-stencil}. Writing
$\lambda^{k+1}$ in its displacement and velocity blocks as before,
\[
\big(\lambda^{k+1}\big)^\top\frac{\partial A}{\partial p_j}x^k
= \Delta t\,\big(\Delta t\,\lambda_u^{k+1} + \lambda_v^{k+1/2}\big)^\top
\frac{\partial \mathbf{L}_{\mathrm{stag}}}{\partial p_j}\,u^k,
\]
and the reduced gradient is assembled, one time step at a time, as
\begin{equation}\label{eq:grad-stag}
(\nabla_p J)_j = \frac{\partial J}{\partial p_j}
- \Delta t\sum\nolimits_{k=0}^{N-1}
\big(\Delta t\,\lambda_u^{k+1} + \lambda_v^{k+1/2}\big)^\top
\frac{\partial \mathbf{L}_{\mathrm{stag}}}{\partial p_j}\,u^k.
\end{equation}
This is the discrete \emph{imaging condition}: each parameter's gradient
is a zero-lag cross-correlation, accumulated over the whole record,
between the forward displacement field $u^k$ and the adjoint field
$\lambda_u^{k+1},\lambda_v^{k+1/2}$ at the same spatial stencil the
forward operator used to couple them---the same object full-waveform-inversion
codes assemble to update a velocity model from a forward and an adjoint
wavefield. The two constructions differ exactly as
Section~\ref{sec:collocated} foreshadowed: the collocated gradient
in~\eqref{eq:grad-collocated} correlates $u^k$ against the single
adjoint field $\lambda_v^{k+1}$ alone, while the staggered gradient
above correlates it against the combination $\Delta
t\,\lambda_u^{k+1}+\lambda_v^{k+1/2}$---one more instance of the
discretization, not the physics, deciding the shape of the computation.

\subsection{Numerical Simulation}\label{sec:numerical-demonstration}
A different discretization does not change the continuous equation
being solved; it changes the structure of the numerical model, and
with it the adjoint. Such choices are typically made for numerical
fidelity---stability limits, accuracy order, or structure-preserving
properties such as energy conservation---rather than for any change in
the underlying physics. 
The collocated scheme updates all variables in a single step, so its
adjoint is a single backward recursion; the staggered scheme splits the
update into two interleaved half-steps, so its adjoint must do the same
in reverse, propagating sensitivities between integer and half-integer
time levels. In both cases the forward computation propagates the state
forward in time while the adjoint propagates sensitivity of the
objective backward through the same steps---the adjoint is determined
by the discrete algorithm, not by the differential equation, and a
change in discretization changes the adjoint equations even though the
underlying physics does not.

\begin{figure}[t]
\centering
\begin{subfigure}[t]{0.48\textwidth}
  \centering
  \includegraphics[width=\linewidth]{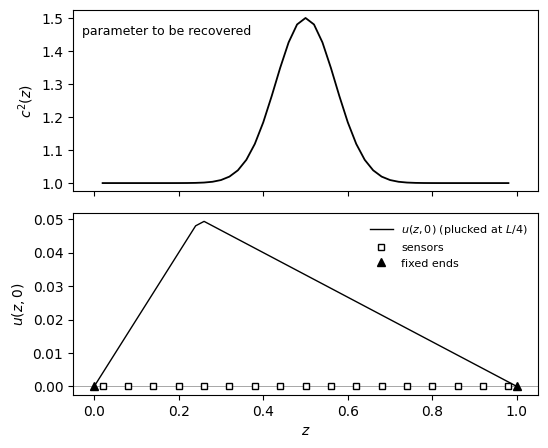}
  \caption{Experiment: the profile $c^2(z)$ to be recovered (top) and the
  plucked string at $t=0$ with the sensor locations (bottom).}
  \label{fig:demo-setup}
\end{subfigure}\hfill
\begin{subfigure}[t]{0.48\textwidth}
  \centering
  \includegraphics[width=\linewidth]{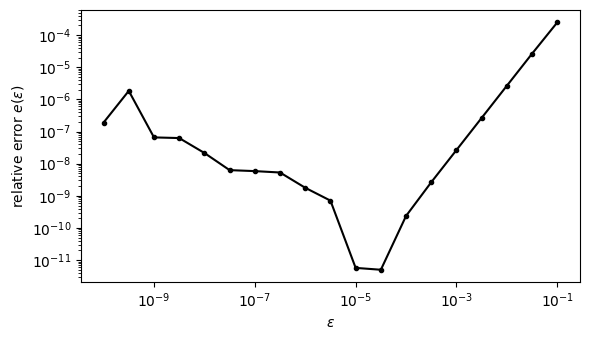}
  \caption{Gradient verification: relative error between the adjoint
  directional derivative and a central finite difference, as a function
  of the step $\varepsilon$.}
  \label{fig:demo-gradcheck}
\end{subfigure}

\medskip
\begin{subfigure}[t]{0.48\textwidth}
  \centering
  \includegraphics[width=\linewidth]{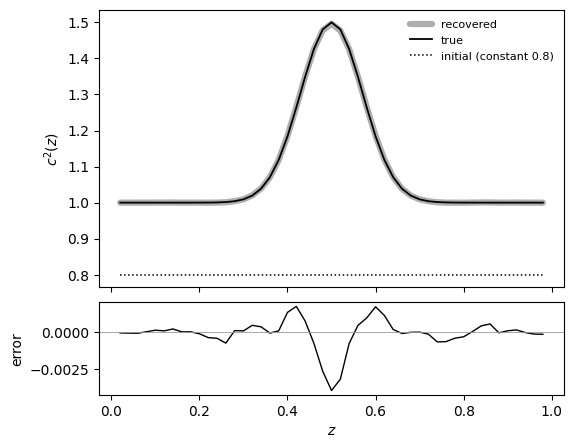}
  \caption{Recovered profile (thick gray) against the true profile
  (black) and the constant initial guess (dotted); the lower panel shows
  the pointwise error.}
  \label{fig:demo-recovery}
\end{subfigure}\hfill
\begin{subfigure}[t]{0.48\textwidth}
  \centering
  \includegraphics[width=\linewidth]{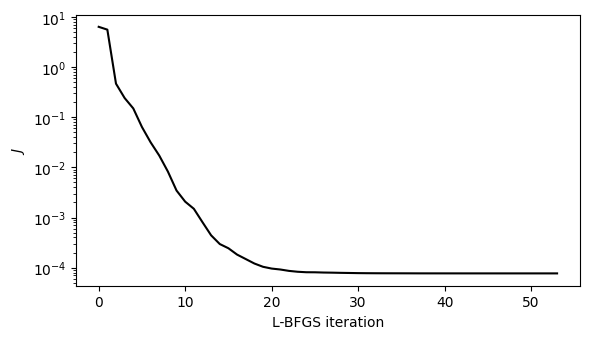}
  \caption{Misfit $J$ against L-BFGS iteration; the plateau is the
  regularization floor $\tfrac{\beta}{2}\mathcal{R}(p_{\mathrm{true}})$.}
  \label{fig:demo-convergence}
\end{subfigure}
\caption{Numerical demonstration with the collocated discretization
($n_z=51$, $n_t\approx6\times10^{3}$), produced by the accompanying
notebook~\cite{Kia2026Colab}. All panels are generated by the code as written from scratch
in this section: the forward sweep $x^{k+1}=A(p)x^k$, the backward
adjoint recursion, the gradient assembly, the finite-difference check,
and the L-BFGS iteration.}
\label{fig:demo}
\end{figure}

The worked example above was kept deliberately small so that every vector and matrix could be
written out by hand. To see the construction at work on a problem large
enough to resemble practice, an accompanying notebook~\cite{Kia2026Colab}
carries out the same steps, with no library solver or optimizer, on a
finer grid: a string of length $L=1$ discretized with $n_z=51$ spatial
nodes. The step $\Delta t$ is set to a small fraction of the CFL limit
$\Delta z/c_{\max}$---not because forward Euler has a stable region
for this undamped system, it does not, but to keep the resulting
growth negligible over the record, as quantified below---and the
record is run long enough, $T=2$, for the disturbance to cross the
domain and reflect off both fixed ends several times; together these
fix $n_t=6123$. Figure~\ref{fig:demo}
reports the results for the collocated discretization, the single case
study carried in the text; the staggered discretization of
Section~\ref{sec:staggered} is built and verified in the same notebook~\cite{Kia2026Colab},
and the closing remark of this section reports what changes, and what
does not, when it replaces the collocated~scheme.

The experiment is shown in Fig.~\ref{fig:demo-setup}. A smooth profile
$c(z)^2$---a localized bump on a unit background, standing in for an
unknown material inhomogeneity---is prescribed as the ``true'' parameter,
with the wave speed at the two fixed ends treated as known boundary data,
exactly as in Section~\ref{sec:collocated}. The string is plucked as in
the worked example at a quarter of its length, so that all modes are excited and the two sides of the bump are
probed differently, and it is released from rest. The forward model then
generates noise-free synthetic displacement records at sensors placed at
every third interior node (seventeen sensors in all). The objective is
the misfit~\eqref{eq:objective} with the regularization term instantiated
as a first-difference smoothness penalty,
\[
\mathcal{R}(p) \;=\; \sum_{j}\left(\frac{p_{j+1}-p_j}{\Delta z}\right)^{2},
\qquad
\frac{\partial J}{\partial p} \;=\; \frac{\beta}{2}\,\nabla\mathcal{R}(p),
\]
with a small weight $\beta=10^{-6}$; the penalty discourages node-to-node
oscillation in the estimate, and it enters the computed gradient exactly
as the $\partial J/\partial p$ term of the reduced-gradient formula, alongside the adjoint contribution $(\partial F/\partial p)^\top\lambda$.
The same discretization is used both to generate this synthetic record
and to run the inversion---a deliberate instance of what is sometimes
called an \emph{inverse crime}\footnote{Generating synthetic data with
the same numerical model used for the inversion can mask errors that
would appear if the two differed; it remains a standard first check of
an inversion algorithm's internal consistency, and is named here to
avoid any implication that a real, independently generated record would
behave identically~\cite{KaipioSomersalo2005}.}, made to isolate the
behavior of the adjoint gradient and the optimizer from the separate
question of model error between a real measurement process and the
model used to invert it.

Before the gradient is used, it is verified, and Fig.~\ref{fig:demo-gradcheck}
is the finite-difference test of Appendix~\ref{sec:gradient_verification}
applied to the code as written. For a fixed direction $\delta p$, the
adjoint prediction $\nabla_pJ^\top\delta p$ is compared with the central
difference $[J(p+\varepsilon\delta p)-J(p-\varepsilon\delta p)]/2\varepsilon$,
which uses only forward solves. The characteristic V-shape is the
verdict: on the right the error decreases with slope two as $\varepsilon$
shrinks, the truncation error of a central difference; on the left it
grows again as round-off in the difference of two nearly equal objective
values is amplified by the division by $\varepsilon$; and at the bottom,
near $\varepsilon\approx10^{-5}$, the two agree to some twelve digits. An
error in the adjoint---a missing transpose, an off-by-one in the backward
recursion, an omitted term in $(\partial F/\partial p)^\top\lambda$---would
flatten this curve at the level of the error, with no descent on the
right, and the test would fail regardless of $\varepsilon$.

The verified gradient drives an L-BFGS iteration, implemented in the
notebook from the two-loop recursion of Section~\ref{sec:solver}.
Section~\ref{sec:solver}'s algorithm states the step length in terms of
a line search satisfying the Wolfe conditions, for generality; the
notebook instead uses a simpler backtracking search enforcing only the
sufficient-decrease (Armijo) condition, and it is this simpler search
that produces the results reported below. The iteration also carries a
projection onto $c^2\geq0.05$ that keeps the wave speed physical:
because a hard projection can produce a curvature pair with $s^\top
y\le0$ on the rare iterations where the bound is active, violating the
assumption the two-loop recursion relies on and risking a non-descent
direction, the notebook checks $s^\top y$ before every update and, when
it is not safely positive, discards the pair and resets the L-BFGS
memory rather than accept it; a bound-constrained variant such as
L-BFGS-B~\cite{Byrd1995} handles this more systematically in settings
where the projection is active often enough to matter. Two standard
stopping rules end the iteration: it terminates when the gradient norm
falls below $10^{-6}$ of its initial value or when $J$ changes by less
than one part in $10^{9}$ over three consecutive iterations.

Two remarks connect this demonstration back to the tutorial's argument.
First, the time step: forward Euler applied to this undamped system
has no stable region at all---every mode is amplified by a factor
$\sqrt{1+\omega^2\Delta t^2}>1$ each step, regardless of $\Delta t$---so
stability is not on offer here, only a growth rate kept small enough to
be negligible over the finite record simulated. Forward Euler is used
in this section for its algebraic simplicity, not because it is an
integrator one would choose for a wave equation in practice; an
implicit scheme such as the trapezoidal (Crank--Nicolson) rule is
unconditionally stable for this system and would be the appropriate
choice in production use, at the cost of a tridiagonal solve per step
whose adjoint is itself a transposed tridiagonal solve. For the
parameters used here ($n_z=51$, $\Delta t\approx3.3\times10^{-4}$,
$n_t=6123$), the growth is real but concentrated where it does the
least damage: the plucked initial condition is smooth, and its three
lowest spatial modes, carrying $99.6\%$ of its initial energy, grow by
at most $3\%$ over the full record; only the poorly resolved, high-
frequency end of the spectrum grows substantially---the highest
resolvable mode by a factor of roughly $10^2$---and that end carries a
negligible share of the energy this excitation contains in the first
place. This is a property of the discrete algorithm, not of the
string, and the staggered scheme does not share it---its leapfrog
update is stable and non-dissipative at the Courant--Friedrichs--Lewy
scale itself, and the notebook confirms this directly: run at matched
time-step and sensor resolution, the staggered scheme's adjoint
gradient---built from the structurally different recursion of
Section~\ref{sec:staggered}---passes the same finite-difference test
and drives the same L-BFGS iteration to a profile of comparable
accuracy. Taking advantage of the larger step the staggered scheme can
stably tolerate is not, however, free: pushed to a step an order of
magnitude larger with the same sensor budget, the notebook's recovery
degrades measurably, since the coarser record carries less information
about the wavefield. The trade is real and worth making
deliberately---fewer time steps at a modest cost in accuracy, or the
reverse---rather than assumed away.

Second, scope: the outcome of an inversion depends on much that this
tutorial deliberately sets aside---the stability and dispersion of the
time integrator, the information content of the excitation and the record
length, the regularization weight, and the convexity of the misfit from
the chosen starting point---and each of these can defeat a correct
gradient. The adjoint's role is narrower and more reliable than any of
them: it delivers the exact gradient of the discretized objective, and
the finite-difference test confirms that it has. When an inversion
disappoints, that test is what allows the adjoint to be ruled out as the
cause, and the practitioner's attention directed to the experiment and
the numerics instead. The notebook is written so that each of these
ingredients can be changed independently, and the reader is encouraged to
do so; it is available at~\cite{Kia2026Colab}.

\clearpage

%%%%%%%%%%%%%%%%%%%
%%%%%%%%%%%%%%%%%%%
\section{Summary and Concluding Remarks}\label{sec:interpretation}

We close this tutorial with a brief recap of the construction developed
above, a discussion of what it leaves out, and a closing perspective.

\noindent\textbf{Summary.}
Section~\ref{sec:inverse} cast the inverse problem as constrained
optimization: the parameters $p$ are free design variables, the state
$x$ is determined from them by the governing equations $F(x,p)=0$, and
after discretization $x$ collects the state over every time level into
a single vector. Section~\ref{sec:adjoint} showed that the reduced
gradient
\[
\nabla_p J = \frac{\partial J}{\partial p} +
\left(\frac{\partial F}{\partial p}\right)^{\!\top}\!\lambda,
\qquad
\left(\frac{\partial F}{\partial x}\right)^{\!\top}\!\lambda
= -\frac{\partial J}{\partial x},
\]
requires only one auxiliary solve, the adjoint equation, in place of the
$m$ forward solves a naive finite-difference gradient would need---the
cost of the gradient is essentially independent of the number of
parameters. Section~\ref{sec:backward} showed why that adjoint solve
runs backward in time: the same block lower-triangular structure that
lets the forward problem be marched forward in time makes its
transpose block upper-triangular, so $\lambda$, matching $x$ in
dimension at every step, is found by marching backward from a terminal
condition fixed by the objective. Section~\ref{sec:solver} showed how
this gradient is consumed by a gradient-based optimizer, and
Section~\ref{sec:discrete-wave} carried the whole construction through
in full detail for two distinct discretizations of the same physical
system, confirming that the adjoint is determined by the discrete
algorithm rather than by the differential equation it approximates.

\vspace{0.5em}
\noindent
\textbf{Scope and extensions.}
The tutorial and its numerical demonstration treat data assimilation as a
deterministic least-squares problem: the objective~\eqref{eq:objective}
penalizes the mismatch between simulated and observed data with no
explicit account of measurement noise, and the parameter is reported as a
single point estimate $p^\star$ rather than as a distribution. Two
consequences of this choice are worth stating plainly, since they bound
what the method presented here can and cannot deliver.

First, real sensor data are noisy, and the least-squares
objective~\eqref{eq:objective} is the negative log-likelihood of the data
under an implicit assumption of independent Gaussian measurement noise.
When that noise is significant, a point estimate obtained by minimizing
it, with or without the Tikhonov-type regularization used
here~\cite{TikhonovArsenin1977,EnglHankeNeubauer1996}, does not by itself
convey how much confidence to place in the recovered parameter. The
Bayesian formulation of the inverse problem addresses this directly: the
parameter is modeled as a random variable, the noise model enters through
a likelihood, and the solution is the full posterior distribution rather
than a single minimizer~\cite{Tarantola2005,KaipioSomersalo2005,Stuart2010}.
The adjoint gradient developed in this tutorial does not disappear in
that setting; it is the same object that drives the maximum a posteriori
estimate and that makes gradient-based sampling and variational posterior
approximation tractable for the large, PDE-constrained parameter spaces
typical of this class of problem~\cite{BuiThanhEtAl2013}.

Second, the demonstration places a fixed, comparatively dense array of
sensors along the string. When measurements are few, or concentrated away
from the region a parameter change most affects, the data no longer
constrain every direction in parameter space: the map from parameters to
data develops a null space, or a near-null space corrupted by noise, and
the inverse problem becomes non-unique or unstable in the sense of
Hadamard~\cite{EnglHankeNeubauer1996}. The smoothness
penalty~$\mathcal{R}(p)$ used in Section~\ref{sec:numerical-demonstration}
is one remedy, chosen for its simplicity; it resolves the
non-uniqueness by fiat, imposing a preference for smooth profiles that
may or may not reflect the true one. A Bayesian prior serves the same
regularizing role in a way that is explicit about the assumption being
made, propagates the resulting uncertainty into the posterior rather
than hiding it in a point estimate, and can encode structure, such as
sparsity or spatial correlation, well beyond what a first-difference
penalty allows~\cite{Stuart2010,KaipioSomersalo2005}.

Neither extension changes the role of the adjoint method developed here;
it changes what the gradient is used for; once instead of iteratively
inside a sampler, or applied to a random realization of the noise instead
of a single data set. Readers pursuing quantified uncertainty or
sparse-data regimes are encouraged to treat this tutorial as the
deterministic core on which that machinery is built.

\vspace{0.5em}
\noindent
\textbf{Closing perspective.}
The adjoint method is a structured way of reusing the governing physics
of a system to compute sensitivities efficiently. Its strength lies not
in the formulas themselves, but in how it transforms a high-dimensional
gradient computation into a pair of coordinated solves: a forward problem
that propagates the state, and an adjoint problem that propagates
sensitivity information backward. Once the roles of the forward model, the adjoint variables, and the
gradient construction are clearly understood, the method becomes
considerably more transparent. Mastering it is less about memorizing
derivations and more about recognizing structure: how information flows
through the system, how dependencies are enforced, and how the same
operators that govern the physics can be reused for sensitivity analysis.
This perspective is what makes the adjoint method practical in
large-scale computational problems, and it is what makes it a way of
thinking about how information flows through a system, not merely a
derivation to be reproduced.

\bibliographystyle{IEEEtran}
\bibliography{references}

\appendix

\section{Verifying the Computed Gradient}
\label{sec:gradient_verification}
% ---------------------------------------------------------------------
 
Equation~\eqref{eq:reduced_grad2} is short, but the code that
implements it is not. Between the formula and a working gradient lie
the assembly of $\partial F/\partial x$ and $\partial F/\partial p$ for a
particular discretization, the treatment of boundary and initial
conditions, the ordering of the state vector, the direction of the time
recursion, and the transpose of every operator involved. Each of these
is an opportunity for a sign error, a dropped transpose, an off-by-one
time index, or a misplaced boundary term. The purpose of this
subsection is to describe how a computed gradient is verified, why the
verification cannot be skipped, and what it reveals when it fails.
 
\paragraph{Why verification is necessary}
An error in the gradient does not announce itself. The optimization
algorithm does not crash; it continues to run, and in most cases the
objective continues to decrease. A gradient with a sign error in one
block, a missing regularization term, or a transposed operator is still
frequently a descent direction, so a line search will still find steps
that reduce $J$. What the practitioner observes is not failure but
degraded performance: slow convergence, a stall at an objective value
that is not a minimum, a reconstruction that fits the data poorly, or
a quasi-Newton method whose curvature pairs are inconsistent and whose
search directions deteriorate. These symptoms are indistinguishable from
those of an ill-conditioned problem, an inadequate regularization, or a
poor initial guess. Without an independent check, the gradient is the
one component of the pipeline whose correctness is assumed rather than
established.
 
Two further features of the adjoint method make verification
particularly important. First, the forward solver and the adjoint
solver are separate pieces of code that must be exact transposes of one
another; nothing in the forward code enforces this, and the two codes
are often written at different times. Second, as emphasized throughout
this tutorial, the adjoint is determined by the discretized algorithm.
An adjoint derived for one discretization and applied to another, or
an adjoint derived from the continuous equations and then discretized
independently, produces a gradient that is \emph{approximately} but not
exactly the gradient of the discrete objective that the code evaluates.
The tests below distinguish an exact discrete gradient from an
approximately consistent one, and this distinction is precisely the
point of working at the discrete level.
 
Verification is inexpensive. The tests that follow require a handful of
forward solves and one adjoint solve, performed once during development.
Against the cost of an optimization run that fails silently, this is
negligible.
 
% ---------------------------------------------------------------------
\subsection{Finite-Difference Test of the Directional Derivative}
% ---------------------------------------------------------------------
 
The most direct check compares the adjoint gradient against a
finite-difference approximation of the derivative of the reduced
objective $J(p)$. Comparing every component $\partial J/\partial p_j$
individually requires $2m$ forward solves, which is affordable for the
small examples of this tutorial but not in general. The check is
therefore carried out along a direction. Let $\delta p\in\mathbb R^m$ be
a fixed vector, chosen at random with entries of order one, and let
$g=\nabla_p J(p)$ denote the adjoint gradient at the current
parameter $p$. The directional derivative of $J$ along $\delta p$ is
$g^{\top}\delta p$, and its central finite-difference approximation is
\begin{equation}\label{eq:fd_directional}
  D_\varepsilon(p;\delta p)
  \;=\;
  \frac{J(p+\varepsilon\,\delta p)-J(p-\varepsilon\,\delta p)}{2\varepsilon},
\end{equation}
which costs two forward solves. The quantity reported is the relative
error
\begin{equation}\label{eq:fd_relerr}
  e(\varepsilon)
  \;=\;
  \frac{\bigl|D_\varepsilon(p;\delta p)-g^{\top}\delta p\bigr|}
       {\bigl|g^{\top}\delta p\bigr|}.
\end{equation}
 
The step $\varepsilon$ should not be fixed at a single value. The error
in~\eqref{eq:fd_directional} has two sources: the truncation error of
the central difference, which behaves as $\mathcal O(\varepsilon^2)$,
and the round-off error from subtracting two nearly equal objective
values, which behaves as $\mathcal O(u/\varepsilon)$ where $u\approx
10^{-16}$ is the double-precision unit round-off. Their sum is minimized
near $\varepsilon\approx u^{1/3}\approx 10^{-5}$, and the smallest
attainable relative error is of order $u^{2/3}\approx 10^{-10}$. Sweeping
$\varepsilon$ over several decades, for instance
$\varepsilon\in\{10^{-1},10^{-2},\dots,10^{-9}\}$, and plotting
$e(\varepsilon)$ on logarithmic axes produces a characteristic
V-shaped curve: the error decreases with slope~$2$, reaches a minimum of
roughly $10^{-8}$ to $10^{-10}$, and then increases as round-off
dominates. This shape is itself the evidence of correctness. A correct
gradient produces the V; an incorrect gradient produces a curve that
flattens at a level of order one and never descends, because the
discrepancy between $D_\varepsilon$ and $g^{\top}\delta p$ is then a
genuine difference, not a discretization artifact.
 
A closely related and more robust form of the same test examines the
Taylor remainder directly. Since
$J(p+\varepsilon\delta p)=J(p)+\varepsilon\,g^{\top}\delta p+
\mathcal O(\varepsilon^2)$ when $g$ is the exact gradient, the residual
\begin{equation}\label{eq:taylor_remainder}
  r(\varepsilon)
  \;=\;
  \bigl|\,J(p+\varepsilon\,\delta p)-J(p)-\varepsilon\,g^{\top}\delta p\,\bigr|
\end{equation}
must decrease as $\varepsilon^2$: on logarithmic axes, successive
halvings of $\varepsilon$ reduce $r$ by a factor of four, and the
slope of $\log r$ against $\log\varepsilon$ is~$2$. If $g$ is wrong in
any component that $\delta p$ excites, the first-order term does not
cancel and $r(\varepsilon)$ decreases only linearly, with slope~$1$.
The test therefore reduces to reading a slope, and a slope of~$1$ where
a slope of~$2$ is expected is an unambiguous indication of an error in
the gradient. Because $\delta p$ is random, a single direction exercises
every component of $g$ simultaneously; two or three independent
directions suffice in practice.
 
Both tests must be run at a parameter value where the gradient is not
small. At or near a minimizer, $g^{\top}\delta p\approx 0$ and the
relative error~\eqref{eq:fd_relerr} is meaningless. A randomly perturbed
parameter, or the initial guess of the optimization, is appropriate.
The regularization weight $\beta$ should be nonzero so that the
$\partial J/\partial p$ term is exercised, and the test should be
repeated at more than one $p$, since a term that vanishes at a
particular parameter value (for instance, a boundary contribution when
the boundary parameters happen to be zero) is not tested there.
 
% ---------------------------------------------------------------------
\subsection{The Adjoint (Dot-Product) Test}
% ---------------------------------------------------------------------
 
The finite-difference test verifies the gradient as a whole but does
not indicate where an error lies. The adjoint test isolates the
component that is most specific to the method: whether the backward
recursion is in fact the transpose of the forward recursion.
 
For any matrix $M$ and vectors $v$, $w$ of compatible dimensions, the
identity
\begin{equation}\label{eq:dot_product_test}
  \langle Mv,\,w\rangle \;=\; \langle v,\,M^{\top}w\rangle
\end{equation}
holds exactly. In an adjoint code the two sides of this identity are
computed by two different routines: the left by the forward
(linearized) solver applied to $v$, the right by the adjoint solver
applied to $w$. If the two routines are correct transposes of one
another, the two inner products agree to round-off, that is, to a
relative difference of order $10^{-14}$ or smaller. The test requires
no step size, involves no truncation error, and passes or fails
decisively.
 
In the setting of this tutorial the test applies at two levels. At the
level of a single time step, $M$ is the one-step operator $A(p)$ of
Section~\ref{sec:discrete-wave} and the test compares the routine that
advances the state with the routine that propagates the adjoint
backward one step; the test is performed with $v$ and $w$ random and
takes microseconds. At the level of the full trajectory, $M$ is the
block matrix $\partial F/\partial x$, the left side of
\eqref{eq:dot_product_test} is computed by running the forward
time-stepping code on a random state perturbation $v$, and the right
side by running the backward sweep on a random $w$. This is the test
that catches an off-by-one time index, a boundary node handled
differently in the two directions, or, in the staggered discretization,
a half-step update applied in the wrong order. Because
$\partial F/\partial x$ is block lower-triangular, an error in a single
block appears in the test as a discrepancy of order one; the identity
is not approximately satisfied when the transpose is approximately
correct.
 
The same test applied to $\partial F/\partial p$ and its transpose
verifies the second term of~\eqref{eq:reduced_grad2}, and applied
to a symmetric operator such as the spatial stencil $\mathbf L(c)$ it
confirms the symmetry that permits $\mathbf L^{\top}$ to be replaced by
$\mathbf L$ in the adjoint recursion.
 
% ---------------------------------------------------------------------
\subsection{A Verification Protocol}
% ---------------------------------------------------------------------
 
The two tests are complementary. The adjoint test is exact and local:
it verifies the transpose relationship between specific pairs of
routines but says nothing about whether $\partial F/\partial x$ itself
has been derived correctly, since a consistent pair of wrong operators
passes it. The finite-difference test is approximate and global: it
verifies the entire gradient against the objective the code actually
evaluates but does not localize an error. Used together, they bound the
problem: an adjoint code that passes both has a forward linearization,
a transpose, and a gradient assembly that are mutually consistent and
consistent with the objective.
 
\begin{tcolorbox}[title=Gradient Verification Protocol,
  colback=gray!5, colframe=black!70, sharp corners]
Given a parameter $p$ (not near a minimizer), regularization
$\beta>0$, and random directions $\delta p\in\mathbb R^m$,
$v,w\in\mathbb R^n$:
\begin{enumerate}
  \item \textbf{Symmetry of the stencil.} Verify
        $\langle\mathbf L v_1, v_2\rangle=\langle v_1,\mathbf L v_2\rangle$
        to round-off.
  \item \textbf{One-step adjoint test.} Verify
        $\langle A v,\,w\rangle=\langle v,\,A^{\top}w\rangle$
        using the forward-step and backward-step routines.
  \item \textbf{Full-trajectory adjoint test.} Verify
        $\bigl\langle \tfrac{\partial F}{\partial x}v,\,w\bigr\rangle
        =\bigl\langle v,\,\tfrac{\partial F}{\partial x}^{\!\top}w\bigr\rangle$
        using the forward sweep and the backward sweep.
  \item \textbf{Parameter Jacobian test.} Verify
        $\bigl\langle \tfrac{\partial F}{\partial p}\delta p,\,w\bigr\rangle
        =\bigl\langle \delta p,\,\tfrac{\partial F}{\partial p}^{\!\top}w\bigr\rangle$.
  \item \textbf{Taylor-remainder test.} Compute $r(\varepsilon)$
        in~\eqref{eq:taylor_remainder} for
        $\varepsilon=10^{-1},\dots,10^{-8}$ and confirm slope~$2$ on
        logarithmic axes.
  \item \textbf{Central-difference test.} Compute $e(\varepsilon)$
        in~\eqref{eq:fd_relerr} over the same range, confirm the
        V-shaped curve, and record the minimum (expected
        $10^{-8}$ to $10^{-10}$).
  \item Repeat steps 5--6 for a second $\delta p$ and a second $p$.
\end{enumerate}
Steps 1--4 pass or fail at round-off. Steps 5--6 pass when the
predicted asymptotic behavior is observed.
\end{tcolorbox}
 
Table~\ref{tab:verification_failures} lists the errors most frequently
encountered when implementing the adjoint of a time-stepping scheme and
the signature each leaves in the tests. The value of the table lies in
the fact that these signatures are distinct: the pattern of which tests
pass and which fail, together with the magnitude of the failure,
usually identifies the error without further debugging.
 
\begin{table}[t]
\centering
\small
\caption{Common implementation errors and their signatures in the
verification tests. ``Plateau'' denotes a relative error that does not
decrease with $\varepsilon$; ``slope~1'' refers to the
Taylor-remainder test~\eqref{eq:taylor_remainder}.}
\label{tab:verification_failures}
\begin{tabular}{@{}p{0.34\linewidth}p{0.30\linewidth}p{0.28\linewidth}@{}}
\toprule
\textbf{Error} & \textbf{Adjoint tests (1--4)} & \textbf{FD tests (5--6)} \\
\midrule
Sign error in the adjoint right-hand side
  & Pass
  & $e(\varepsilon)\approx 2$; slope~1 \\
Transpose omitted on a non-symmetric block
  & Fail at step 2 or 3, $\mathcal O(1)$ discrepancy
  & Plateau; slope~1 \\
Off-by-one time index in the backward sweep
  & Pass at step 2, fail at step 3
  & Plateau; slope~1 \\
Boundary node treated differently forward and backward
  & Fail at step 3, discrepancy $\propto$ boundary contribution
  & Plateau at a level that depends on $p$ \\
$\partial J/\partial p$ (regularization) term omitted
  & Pass
  & Plateau at a level $\propto\beta$; vanishes when $\beta=0$ \\
Staggered half-step updates applied in the wrong order
  & Fail at step 2
  & Plateau; slope~1 \\
Adjoint of a different discretization used
  & Pass (self-consistent) or fail, depending on implementation
  & Plateau at $\mathcal O(\Delta t)$ or $\mathcal O(\Delta z)$; decreases only under mesh refinement \\
Initial state $x^0$ given an adjoint variable and included in the gradient
  & Pass
  & Plateau; slope~1 \\
Correct gradient, $\varepsilon$ too large or too small
  & Pass
  & V-shape present, minimum not reached at the tested $\varepsilon$ \\
\bottomrule
\end{tabular}
\end{table}
 
The penultimate row deserves comment. When the adjoint is derived for a
discretization other than the one implemented in the forward solver,
or is obtained by discretizing the continuous adjoint equation
independently, the resulting gradient is consistent with the discrete
objective only up to the discretization error. The finite-difference
test then plateaus at a level of order $\Delta t$ or $\Delta z$ rather
than descending to round-off, and the plateau moves only when the mesh
is refined. This behavior is not a bug in the usual sense, and such
gradients are used deliberately in some settings; but it is invisible
without the test, and it is the reason the gradient produced by the
discrete-then-optimize route of this tutorial passes the test to
$10^{-10}$ while an optimize-then-discretize gradient for the same
problem does not.
 
\subsection{Verification when the adjoint solver is not accessible}
The finite-difference tests require only the ability to evaluate the
objective at a perturbed parameter and to obtain a gradient; neither
requires access to the internals of the adjoint code. They can
therefore be applied to a gradient produced by third-party software,
using that software's own objective evaluation, to establish whether
the gradient it returns is the gradient of the objective it computes.
For the practitioner who uses such software without the ability to
inspect it, this is the one test that connects the reported sensitivity
to the reported objective, and it should be the first thing done with
any new problem setup.
 
\subsection{Summary}
The gradient is the single point of contact between the governing
equations and the optimization algorithm. Every subsequent
result---convergence history, reconstructed parameters, reported
sensitivities---rests on it. The tests of this subsection establish its
correctness at a cost of a few forward solves, they distinguish an
exact discrete gradient from an approximately consistent one, and when
they fail their pattern of failure indicates where to look. They are
applied to both discretizations of Section~\ref{sec:discrete-wave},
where the staggered case in particular illustrates how an error in
half-step bookkeeping that is invisible in the forward solution is
caught immediately by the one-step adjoint test.

\end{document}